\long\def\ignore#1\endignore{}
\long\def\later#1{\red XXX\endred}
\def\blackpen{\pdfliteral{.0 .0 .0 rg .0 .0 .0 RG}}
\def\redpen{\pdfliteral{.9 .0 .0 rg .9 .0 .0 RG}}
\long\def\red#1\endred{\redpen#1\blackpen}
\font\tenmsa=msam10
\font\tenmsb=msbm10

\font\largebf=cmbx10 scaled\magstep2
\def\all{\hbox{for all}}
\def\and{\hbox{and}}
\def\bra#1#2{\langle#1,#2\rangle}
\def\Bra#1#2{\big\langle#1,#2\big\rangle}
\def\cite#1\endcite{[#1]}
\def\dist{\hbox{\rm dist}}
\def\dom{\hbox{\rm dom}}

\def\dbs{^{**}}
\def\episum{\mathop{\nabla}}
\def\eps{\varepsilon}
\def\f#1#2{{#1 \over #2}}
\def\fbar{\overline f}
\def\fourth{\ts\f14}
\def\half{{\textstyle\f12}}

\def\infn{\inf\nolimits}
\def\limn{\lim\nolimits}

\def\lr{\Longrightarrow}

\def\minn{\min\nolimits}
\def\M{{\cal M}}

\def\on{\hbox{on}}
\def\P{{\cal P}}
\def\PC{{\cal PC}}
\def\PCLSC{{\cal PCLSC}}
\def\PCLSCST{{\cal PCLSC}^*}
\def\pt{\wt p}
\def\qed{\hfill\hbox{\tenmsa\char03}}
\def\qlr{\quad\lr\quad}

\def\qt{\wt q}

\def\quand{\quad\and\quad}
\def\r{\hbox{\tenmsb R}}
\def\rbar{\,]{-}\infty,\infty]}

\def\rthalf{\ts\f1{\sqrt2}}
\def\rttwo{\sqrt2}
\def\Sem{{\cal S}}
\long\def\slant#1\endslant{{\sl#1}}

\def\supn{\sup\nolimits}
\def\ts{\textstyle}
\def\wh{\widehat}

\def\wt{\widetilde}
\def\xbra#1#2{\lfloor#1,#2\rfloor}
\def\Xbra#1#2{\big\lfloor#1,#2\big\rfloor}
\def\ybra#1#2{\lceil#1,#2\rceil}

\def\defSection#1{}
\def\defCorollary#1{}
\def\defDefinition#1{}
\def\defExample#1{}
\def\defLemma#1{}
\def\defNotation#1{}
\def\defProblem#1{}
\def\defRemark#1{}
\def\defTheorem#1{}
\def\locno#1{}
\def\meqno#1{\eqno(#1)}
\def\nmbr#1{}
\def\Proof{\medbreak\noindent{\bf Proof.}\enspace}
\def\Proo{{\bf Proof.}\enspace}
\def\Signoff{}
\def \INTsec{0}
\def \SSDsec{1}
\def \PCdef{1.1}
\def \SSDdef{1.2}
\def \Hex{1.3}
\def \EEex{1.4}
\def \NOTSSDex{1.5}
\def \ROOTlem{1.6}
\def \ROOTrem{1.7}
\def \POSFdef{1.8}
\def \Llem{1.9}
\def \FATdef{1.10}
\def \FFATlem{1.11}
\def \FFATrem{1.12}
\def \BSSDsec{2}
\def \NORMdef{2.1}
\def \IOTAone{2.1.1}
\def \IOTAtwo{2.1.2}
\def \IOTAthree{2.1.3}
\def \IOTAfour{2.1.4}
\def \IOTAfive{2.1.5}
\def \IOTAsix{2.1.6}
\def \NORMrem{2.2}
\def \EENex{2.3}
\def \EEone{2.3.1}
\def \EESPECex{2.4}
\def \VZdef{2.5}
\def \VZone{2.5.1}
\def \VZtwo{2.5.2}
\def \VZthree{2.5.3}
\def \Pdef{2.6}
\def \EXlem{2.7}
\def \EXone{2.7.1}
\def \EXtwo{2.7.2}
\def \EXthree{2.7.3}
\def \EXfour{2.7.4}
\def \EXfive{2.7.5}
\def \EXseven{2.7.6}
\def \EXeight{2.7.7}
\def \Alem{2.8}
\def \EXthm{2.9}
\def \MLTrem{2.10}
\def \Hcor{2.11}
\def \THdef{2.12}
\def \THlem{2.13}
\def \THrem{2.14}
\def \THthm{2.15}
\def \THone{2.15.1}
\def \HHATtwo{2.15.2}
\def \HHATthree{2.15.3}
\def \MAXthm{2.16}
\def \FPHIrem{2.17}
\def \NOQTrem{2.18}
\def \BSSDAPPSsec{3}
\def \NTHthm{3.1}
\def \DUALsec{4}
\def \DUALdef{4.1}
\def \DUALone{4.1.1}
\def \DUALtwo{4.1.2}
\def \DUALthree{4.1.3}
\def \DUALfour{4.1.4}
\def \DUALfive{4.1.5}
\def \PSEUDOdef{4.2}
\def \PSEUDOone{4.2.1}
\def \DUALrem{4.3}
\def \EEDUALex{4.4}
\def \RSUMlem{4.5}
\def \SHARPrem{4.6}
\def \SSDDlem{4.7}
\def \Hone{4.7.1}
\def \MASdef{4.8}
\def \EQthm{4.9}
\def \MULTIthm{4.10}
\def \DUALAPPSsec{5}
\def \DIFFrem{5.1}
\def \SPECdef{5.2}
\def \INVARthm{5.3}
\def \INVARdef{5.4}
\def \BRthm{5.5}
\def \NAPAIRone{5.5.1}
\def \TWOrem{5.6}
\def \SRdef{5.7}
\def \SRthm{5.8}
\def \SRone{5.8.1}
\def \FMsec{6}
\def \FMthm{6.1}
\def \FMone{6.1.1}
\def \FMtwo{6.1.2}
\def \FENMORSUFFone{6.1.3}
\def \FENMORSUFFtwo{6.1.4}
\def \FENMORSUFFthree{6.1.5}
\def \FENMORSUFFfour{6.1.6}
\def \BS{1}
\def \FITZ{2}
\def \GOSSEZC{3}
\def \ASTWO{4}
\def \ASTHREE{5}
\def \MT{6}
\def \PENOT{7}
\def \FENCHEL{8}
\def \RANGE{9}
\def \BR{10}
\def \PANDM{11}
\def \HBL{12}
\def \HBM{13}
\def \ZAGRODNY{14}
\def \VZ{15}
\def \ZBOOK{16}
%
\magnification 1200
\headline{\ifnum\folio=1
{\hfil{\largebf Nonreflexive Banach SSD spaces}\hfil}
\else\centerline{\rm {\bf Nonreflexive Banach SSD spaces}}\fi}
\bigskip
\centerline{by S. Simons}
\bigskip
\centerline{\bf Abstract}
\medskip
\noindent
In this paper, we unify the theory of SSD spaces, part of the theory of strongly representable multifunctions, and the theory of the equivalence of various classes of maximally monotone multifunctions.
\defSection \INTsec
\medbreak
\centerline{\bf \INTsec\quad Introduction}
\medskip
\noindent
In this paper, we unify three different lines of investigation:  the theory of SSD spaces as expounded in \cite\PANDM\endcite\ and \cite\HBM\endcite, part of the theory of strongly representable multifunctions as expounded in \cite\VZ\endcite\ and \cite\ASTWO\endcite, and the equivalence of various classes of maximally monotone multifunctions, as expounded in \cite\ASTHREE\endcite.
\par
The purely algebraic concepts of \slant SSD space\endslant\ and \slant$q$--positive set\endslant\ are introduced in Definition \SSDdef.   These were originally defined in \cite\PANDM\endcite, and the development of the theory was continued in \cite\HBM\endcite.   Apart from the fact that we write ``$\P$'' instead of ``pos'', we use the notation of the latter of these references.   We show in Lemma \Llem\ how certain proper convex functions $f$ on an SSD space lead to a \slant $q$--positive\endslant\ set, $\P(f)$. In Definition \FATdef, we define the \slant intrinsic conjugate\endslant, $f^@$, of a proper convex function on an SSD space, and we end Section \SSDsec\ by proving in Lemma \FFATlem\ a simple, but useful, property of intrinsic conjugates. 
\par
In Definition \NORMdef, we introduce the concept of a \slant Banach SSD space\endslant, which is an SSD space with a Banach space structure satisfying the compatibility conditions (\IOTAone) and (\IOTAtwo).   A proper convex function on a Banach SSD space may be a \slant VZ function\endslant, which is introduced in Definition \VZdef.   Our main result on VZ functions, established in Theorem \EXthm(c,d), is that \slant if $f$ is a lower semicontinuous VZ function then $\P(f)$ is {\rm maximally} $q$--positive, $f^@$ is also a VZ function, and $\P\big(f^@\big) = \P(f)$\endslant.   Lemma \EXlem(b) is an important stepping--stone to Theorem \EXthm.   In Definition \THdef\ and Lemma \THlem, we introduce and discuss the properties of various convex functions on a Banach SSD space and its dual, and show in Theorem \THthm(c) that if $f$ is a lower semicontinuous VZ function on a Banach SSD then there is a whole family of VZ functions $h$ associated with $f$ such that $\P(h) = \P(f)$.
\par
If $E$ is a nonzero Banach space then it is shown in Examples \EEex, \EENex, and \EESPECex\ that $E \times E^*$ is a Banach SSD space under various different norms.   We show in Section \BSSDAPPSsec\ how the definitions and results of Section \BSSDsec\ specialize to this case.   Theorem \NTHthm\ extends some concepts and results from \cite\BS\endcite\ and \cite\ASTHREE\endcite.   The definition of VZ function involves the norm of $B$ in an essential way.   Looking ahead, we will see in Theorem \INVARthm\ that there is a large class of norms on $E \times E^*$ for which the classes of VZ functions coincide.   This follows from the analysis in Section \DUALsec, which we will now discuss.     
\par
In Definition \DUALdef, we introduce the concept of a \slant Banach SSD dual space\endslant, which is the dual of a Banach SSD space which has an SSD structure in its own right, satisfying the compatibility conditions (\DUALone) and (\DUALtwo).   In this situation, a proper convex function on the (original) Banach SSD space may be an \slant MAS function\endslant, which is introduced in Definition \MASdef.   The main result here is Theorem \EQthm(c), in which we prove that, under the \slant $\pt$--density\endslant\ condition (\PSEUDOone), a function is an MAS function if, and only if, it is a VZ function.   The main stepping stone to Theorem \EQthm\ is Lemma \SSDDlem, which relies on Rockafellar's formula for the conjugate of the sum of two convex functions.
\par
The subtlety of the analysis outlined in the previous paragraph is that definition of MAS function does not use the norm of $B$ explicitly --- it only uses the knowledge of $B^*$.  In a certain sense, the analysis of Section \BSSDsec\ is \slant isometric\endslant, while the analysis of Section \DUALsec\ is \slant isomorphic\endslant, though it would be a mistake to push this analogy too far because, despite the fact that the definition of MAS function does not use the norm of $B$, the conditions (\DUALtwo) and (\PSEUDOone) referred to above do use the norm very strongly.    
\par
In Section \DUALAPPSsec, we show how the results of Section \DUALsec\ specialize to the $E \times E^*$ case.   In Theorem \BRthm, we show how the \slant negative alignment\endslant\ analysis introduced in \cite\BR, Section 8, pp.\ 274--280\endcite\ and \cite\HBM, Section 42, pp.\ 161--167\endcite\ can be used to obtain, and in some cases strengthen, results from \cite\ASTWO\endcite\ and \cite \VZ\endcite.   In Theorem \SRthm, we generalize some equivalencies from \cite\ASTHREE, Theorem 1.2\endcite.   In particular, we give a proof of the very nice result from \cite\ASTHREE\endcite\ that a maximally monotone multifunction is strongly representable if, and only if, it is of type (NI).
\par
At one point in this paper, we will use the Fenchel--Moreau theorem for a not necessarily Hausdorff locally convex space.   For the convenience of the reader, we give a proof of this result in the Appendix, Section \FMsec.
\par
The author would like to thank Constantin Z\u alinescu for making him aware of the preprints \cite\ASTWO\endcite\ and \cite\VZ\endcite, and Benar Svaiter for making him aware of the preprint \cite\ASTHREE\endcite.   He would also like to thank Constantin Z\u alinescu for some very perceptive comments on an earlier version of this paper. 
\defSection\SSDsec
\medbreak
\centerline{\bf \SSDsec\quad SSD spaces}
\medskip
\noindent
We first introduce the concepts of an SSD space and $q$--positive set.   As pointed out in the introduction, these were introduced in \cite\PANDM\endcite\ and \cite\HBM\endcite.  The first of these references has a detailed discussion of the finite dimensional case.
\defDefinition \PCdef
\medbreak
\noindent
{\bf Definition \PCdef.}\enspace If $X$ is a nonzero vector space and $f\colon\ X \to \rbar$, we write $\dom\,f$ for the set $\big\{x \in X\colon\ f(x) \in \r\big\}$.   $\dom\,f$ is the \slant effective domain \endslant of $f$.   We say that $f$ is \slant proper\endslant\ if $\dom\,f \ne \emptyset$.   We write $\PC(X)$ for the set of all proper convex functions from $X$ into $\rbar$.   If $X$ is a nonzero Banach space, we write $\PCLSC(X)$ for the set 
$$\{f \in \PC(X)\colon\ f\ \hbox{is lower semicontinuous on}\ X\},$$
and $\PCLSCST(X^*)$ for the set 
$$\{f \in \PC(X^*)\colon\ f\ \hbox{is $w(X^*,X)$--lower semicontinuous on}\ X^*\}.$$ \par
\defDefinition \SSDdef
\medbreak
\noindent
{\bf Definition \SSDdef.}\enspace We will say that $\big(B,\xbra\cdot\cdot\big)$ is a \slant symmetrically self--dual space (SSD space)\endslant\ (if there is no risk of confusion, we will say simply ``$B$ is an SSD space'') if $B$ is a nonzero real vector space and $\xbra\cdot\cdot\colon B \times B \to \r$ is a symmetric bilinear form.   We define the quadratic form $q$ on $B$ by $q(b) := \half\xbra{b}{b}$.   Let $A \subset B$.   We say that $A$ is \slant$q$--positive\endslant\ if $A \ne \emptyset$ and
$$b,c \in A \lr q(b - c) \ge 0.$$
We say that $A$ is \slant maximally $q$--positive\endslant\ if $A$ is $q$--positive and $A$ is not properly contained in any other $q$--positive set.   We make the elementary observation that if $b \in B$ and $q(b) \ge 0$ then the linear span $\r b$ of $\{b\}$ is $q$--positive.
\medbreak
We now give some examples of SSD spaces and their associated $q$--positive sets.   
\defExample \Hex 
\medbreak
\noindent
{\bf Example \Hex.}\enspace Let $B$ be a Hilbert space with inner product $(b,c) \mapsto \bra{b}{c}$ and $T\colon B \to B$ be a self--adjoint linear operator.   Then $B$ is an SSD space with $\xbra{b}{c} := \bra{Tb}{c}$, and then $q(b) = \half\bra{Tb}{b}$.   Here are three special cases of this example:\par
(a)\enspace If, for all $b \in B$, $Tb = b$ then $\xbra{b}{c} := \bra{b}{c}$, $q(b) = \half\|b\|^2$ and every subset of $B$ is $q$--positive
\par
(b)\enspace If, for all $b \in B$, $Tb = -b$ then $\xbra{b}{c} := -\bra{b}{c}$, $q(b) = -\half\|b\|^2$ and the $q$--positive sets are the singletons.
\par
(c)\enspace If $B = \r^3$ and $T(b_1,b_2,b_3) = (b_2,b_1,b_3)$ then
$$\Xbra{(b_1,b_2,b_3)}{(c_1,c_2,c_3)} := b_1c_2 + b_2c_1 + b_3c_3,$$
and\quad$q(b_1,b_2,b_3) = b_1b_2 + \half b_3^2$.   Here, If $M$ is any nonempty monotone subset of $\r \times \r$ (in the obvious sense) then $M \times \r$ is a $q$--positive subset of $B$.   The set $\r(1,-1,2)$ is a $q$--positive subset of $B$ which is not contained in a set $M \times \r$ for any monotone subset of $\r \times \r$.   The helix $\big\{(\cos\theta,\sin\theta,\theta)\colon \theta \in \r\big\}$ is a $q$--positive subset of $B$, but if $0 < \lambda < 1$ then the helix $\big\{(\cos\theta,\sin\theta,\lambda\theta)\colon \theta \in \r\big\}$ is not.\par
\defExample \EEex 
\medbreak
\noindent
{\bf Example \EEex.}\enspace Let $E$ be a nonzero Banach space and $B := E \times E^*$.   For all $b = (x,x^*)$ and $c = (y,y^*) \in B$, we set $\xbra{b}{c} := \bra{x}{y^*} + \bra{y}{x^*}$.  Then $B$ is an SSD space and
$$q(b) = \half\big[\bra{x}{x^*} + \bra{x}{x^*}\big] = \bra{x}{x^*}.$$
Consequently, if $b = (x,x^*)\ \and\ c = (y,y^*) \in B$ then
$$\bra{x - y}{x^* - y^*} = q(x - y,x^* - y^*) = q\big((x,x^*) - (y,y^*)\big) = q(b - c).$$   
Thus if $A \subset B$ then $A$ is $q$--positive exactly when $A$ is a nonempty monotone subset of $B$ in the usual sense, and $A$ is maximally $q$--positive exactly when $A$ is a maximally monotone subset of $B$ in the usual sense.   We point out that any finite dimensional SSD space of the form described here must have \slant even\endslant\ dimension.  Thus cases of Example \Hex\ with finite odd dimension cannot be of this form.
\defExample \NOTSSDex
\medbreak
\noindent
{\bf Example \NOTSSDex.}\enspace $\r^3$ is \slant not\endslant\ an SSD space with
$$\Xbra{(b_1,b_2,b_3)}{(c_1,c_2,c_3)} := b_1c_2 + b_2c_3 + b_3c_1.$$
(The bilinear form $\xbra\cdot\cdot$ is not symmetric.)
\defLemma \ROOTlem
\medbreak
\noindent
{\bf Lemma \ROOTlem.}\enspace\slant Let $B$ be an SSD space, $f \in \PC(B)$, $f \ge q$ on $B$ and $b,c \in B$.   Then
$$-q(b - c) \le \Big[\sqrt{(f - q)(b)} + \sqrt{(f - q)(c)}\Big]^2.$$\endslant
\Proof We can and will suppose that $0 \le (f - q)(b)  < \infty$ and $0 \le (f - q)(c)  < \infty$.   Let $\sqrt{(f - q)(b)} < \beta < \infty$ and $\sqrt{(f - q)(c)} < \gamma < \infty$, so that $\beta^2 + q(b) > f(b)$ and $\gamma^2 + q(c) > f(c)$.   Then
$$\eqalign{\beta\gamma + \f{\gamma q(b) + \beta q(c)}{\beta + \gamma}
&= \f\gamma{\beta + \gamma}\big(\beta^2 + q(b)\big) + \f\beta{\beta + \gamma}\big(\gamma^2 + q(c)\big)\cr
&> \f\gamma{\beta + \gamma}f(b) + \f\beta{\beta + \gamma}f(c) \ge
f\bigg(\f{\gamma b + \beta c}{\beta + \gamma}\bigg)\cr
&\ge q\bigg(\f{\gamma b + \beta c}{\beta + \gamma}\bigg)
= \f{\gamma^2q(b) + \gamma\beta\xbra{b}{c} + \beta^2q(c)}{(\beta + \gamma)^2}.}$$
Clearing of fractions, we obtain
$$(\beta + \gamma)^2\beta\gamma + (\beta + \gamma)\big(\gamma q(b) + \beta q(c)\big) >
\gamma^2q(b) + \gamma\beta\xbra{b}{c} + \beta^2q(c),$$
from which $(\beta + \gamma)^2\beta\gamma > -\beta\gamma q(b) + \beta\gamma\xbra{b}{c} - \beta\gamma q(c) = -\beta\gamma q(b - c)$.   If we now divide by $\beta\gamma$, we obtain $(\beta + \gamma)^2 > -q(b - c)$, and the result follows by letting $\beta \to \sqrt{(f - q)(b)}$ and $\gamma \to \sqrt{(f - q)(c)}$.\qed 
\defRemark \ROOTrem
\medbreak
\noindent
{\bf Remark \ROOTrem.}\enspace It follows from Lemma \ROOTlem\ and the Cauchy--Schwarz inequality that
$$-q(b - c) \le 2(f - q)(b) + 2(f - q)(c).$$
In the situation of Example \EEex, we recover \cite\VZ, Proposition 1\endcite.
\defDefinition \POSFdef
\medbreak
\noindent
{\bf Definition \POSFdef.}\enspace If $B$ be an SSD space, $f \in \PC(B)$ and $f \ge q$ on $B$, we write
$$\P(f) := \big\{b \in B\colon\ f(b) = q(b)\big\}.$$\par
The following result is suggested by Burachik--Svaiter, \cite\BS, Theorem 3.1, pp. 2381--2382\endcite\ and Penot, \cite\PENOT, Proposition 4\big(h)$\lr$(a), pp. 860--861\endcite.
\defLemma \Llem
\medbreak
\noindent
{\bf Lemma \Llem.}\enspace\slant Let $B$ be an SSD space, $f \in \PC(B)$, $f \ge q$ on $B$ and $\P(f) \ne \emptyset$.   Then $\P(f)$ is a $q$--positive subset of $B$.\endslant
\Proof This is immediate from Lemma \ROOTlem.\qed
\medbreak
We now introduce a concept of conjugate that is intrinsic to an SSD space without any topological conditions.
\defDefinition \FATdef
\medbreak
\noindent
{\bf Definition \FATdef.}\enspace If $B$ is an SSD space and $f \in \PC(B)$, we write $f^@$ for the Fenchel conjugate of $f$ with respect to the pairing $\xbra\cdot\cdot$, that is to say,
$$\all\ c \in B,\qquad f^@(c) := \supn_{b \in B}\big[\xbra{b}{c} - f(b)\big].$$
\par
Our next result represents an improvement of the result proved in \cite\HBM, Lemma 19.12, p.\ 82\endcite, and uses a disguised differentiability argument.   See Remark \FFATrem\ below for another proof of Lemma \FFATlem, due to Constantin Z\u alinescu.   
\defLemma \FFATlem
\medbreak
\noindent
{\bf Lemma \FFATlem.}\enspace\slant Let $B$ be an SSD space, $f \in \PC(B)$ and  $f \ge q$ on $B$.   Then:
$$\leqalignno{a \in \P(f)\ \and\ b \in B &\qlr \xbra{b}{a} \le q(a) + f(b).&(a)\cr
a \in \P(f) &\qlr f^@(a) = q(a).&(b)}$$\endslant
\Proo Let $a \in \P(f)$ and $b \in B$.   Let $\lambda \in \,]0,1[\,$.   For simplicity in writing, let $\mu := 1 - \lambda \in \,]0,1[\,$.   Then
$$\eqalign{\lambda^2q(b) + \lambda\mu\xbra{b}{a} + \mu^2q(a) &= q\big(\lambda b + \mu a\big) \le f(\lambda b + \mu a)\cr
&\le \lambda f(b) + \mu f(a) = \lambda f(b) + \mu q(a).}$$
Thus\quad $\lambda^2q(b) + \lambda\mu\xbra{b}{a} \le \lambda f(b) + \lambda\mu q(a)$.\quad We now obtain (a) by dividing by $\lambda$ and letting $\lambda \to 0$.   Now let $a \in \P(f)$.   From (a), \quad$b \in B \lr \xbra{a}{b} - f(b) \le q(a)$,\quad and it follows by taking the supremum over $b \in B$ that $f^@(a) \le q(a)$.   On the other hand,\quad $f^@(a) \ge \xbra{a}{a} - f(a) = 2q(a) - q(a) = q(a)$,\quad completing the proof of (b).\qed
\defRemark \FFATrem
\medbreak
\noindent
{\bf Remark \FFATrem.}\enspace The author is grateful to Constantin Z\u alinescu for pointing out to him the following alternative proof of Lemma \FFATlem(a).   From Lemma \ROOTlem, with $c$ replaced by $a$,
$-q(b) + \xbra{b}{a} - q(a) = -q(b - a) \le (f - q)(b)$.   Thus $\xbra{b}{a} - q(a) \le f(b)$, as required. 
\defSection\BSSDsec
\medbreak
\centerline{\bf \BSSDsec\quad Banach SSD spaces}
\defDefinition \NORMdef
\medbreak
\noindent
{\bf Definition \NORMdef.}\enspace We say that $B$ is a \slant Banach SSD space\endslant\ if $B$ is an SSD space and $\|\cdot\|$ is a norm on $B$ with respect to which $B$ is a Banach space with norm--dual $B^*$,
$$\half\|\cdot\|^2 + q \ge 0\ \on\ B\meqno\IOTAone$$
and there exists $\iota \in L(B,B^*)$ such that
$$\all\ b,c\in B,\quad \Bra{b}{\iota(c)} = \xbra{b}{c}\hbox{,\quad \big(from which } \big|\xbra{b}{c}\big| \le \|\iota\|\|b\|\|c\|\hbox{\big).}\meqno\IOTAtwo$$
Then, for all $d,e \in B$,
$$|q(d) - q(e)| = \half\big|\xbra{d}{d} - \xbra{e}{e}\big| = \half\big|\xbra{d- e}{d + e}\big| \le \half \|\iota\|\|d - e\|\|d + e\|.\meqno\IOTAthree$$
We define the continuous convex functions $g$ and $p$ on $B$ by $g := \half\|\cdot\|^2$ and $p := g + q$, so that\quad $p \ge 0$ on $B$.\quad Since $p(0) = 0$, in fact
$$\infn_Bp = 0.\meqno\IOTAfour$$
Also, for all $d,e \in B$, $|g(d) - g(e)| = \half\big|\|d\| - \|e\|\big|\big(\|d\| + \|e\|\big) \le \half\|d - e\|\big(\|d\| + \|e\|\big)$.   Combining this with (\IOTAthree), for all $d,e \in B$,
$$|p(d) - p(e)| \le \half\big(1 + \|\iota\|\big)\|d - e\|\big(\|d\| + \|e\|\big).\meqno\IOTAfive$$
(\IOTAtwo) implies that, for all $f \in \PC(B)$ and $c \in B$, $f^@(c) = \supn_{b \in B}\big[\Bra{b}{\iota(c)} - f(b)\big] = f^*\big(\iota(c)\big)$, that is to say,
$$f^@ = f^*\circ\iota\ \on\ B.\meqno\IOTAsix$$\par
\defRemark \NORMrem
\medbreak
\noindent
{\bf Remark \NORMrem.}\enspace  Example \Hex\ is a Banach SSD space provided that $\|T\| \le 1$.   This is the case with (a), (b) and (c) of Example \Hex.
\defExample \EENex
\medbreak
\noindent
{\bf Example \EENex.}\enspace We now continue our discussion of Example \EEex.   We suppose that $B = E \times E^*$ and $\big(B,\|\cdot\|\big)$ is a Banach SSD space such that $B^* = E\dbs \times E^*$, under the pairing
$$\bra{b}{c^*} := \bra{x}{y^*} + \bra{x^*}{y\dbs}\quad\big(b = (x,x^*) \in B,\ c^* = (y\dbs,y^*) \in B^*\big).\meqno\EEone$$
We recall that, for all $(x,x^*) \in B$, $q(x,x^*) = \bra{x}{x^*}$.   It is clear that, for all $(x,x^*) \in B$, $\iota(x,x^*) := (\wh{x},x^*)$ where $\wh{x}$ is the canonical image of $x$ in $E\dbs$.   We note that if $\big(B,\|\cdot\|\big)$ is a Banach SSD space and $\|\cdot\|'$ is a larger norm on $B$ such that $\big(B,\|\cdot\|'\big)^* = E\dbs \times E^*$ then $\big(B,\|\cdot\|'\big)$ is also a Banach SSD space.
\defExample \EESPECex
\medbreak
\noindent
{\bf Example \EESPECex.}\enspace We now discuss some specific examples of the above concepts.   Here it is convenient to introduce a parameter $\tau > 0$.   ($\tau$ stands for ``torsion''.) Then\break $E \times E^*$ is a Banach SSD space if we use the norm $\|(x,x^*)\|_{1,\tau} := \rthalf\big(\tau\|x\| +\|x^*\|/\tau\big)$ or $\|(x,x^*)\|_{2,\tau} := \sqrt{\tau^2\|x\|^2 + \|x^*\|^2/\tau^2}$ or $\|(x,x^*)\|_{\infty,\tau} := \rttwo\big(\tau\|x\| \vee \|x^*\|/\tau\big)$.   (These are arranged in order of increasing size.)   Then the dual norm of $\big(B,\|\cdot\|_{1,\tau}\big)$ is given by $\|(y\dbs,y^*)\|_{\infty,\tau} := \rttwo\big(\tau\|y\dbs\| \vee \|y^*\|/\tau\big)$, the dual norm of $\big(B,\|\cdot\|_{2,\tau}\big)$ is given by  $\|(y\dbs,y^*)\|_{2,\tau} := \sqrt{\tau^2\|y\dbs\|^2 + \|y^*\|^2/\tau^2}$, and the dual norm of $\big(B,\|\cdot\|_{\infty,\tau}\big)$ is given by $\|(y\dbs,y^*)\|_{1,\tau} := \rthalf\big(\tau\|y\dbs\| +\|y^*\|/\tau\big)$.     
\defDefinition \VZdef
\medbreak
\noindent
{\bf Definition \VZdef.}\enspace
Let $X$ be a vector space and $h,k\colon X \to \rbar$.   The inf--convolution of $h$ and $k$ is defined by\quad $(h \episum k)(x) := \infn_{y \in X} \big[h(y) + k(x - y)\big]$\quad($x \in X$).\quad   It is clear that
$$\infn_Xk = 0 \qlr \infn_X\big[h \episum k\big] = \infn_Xh.\meqno\VZone$$
Now let $\big(B,\|\cdot\|\big)$ be a Banach SSD space and $f \in \PC(B)$.   We say that $f$ is a \slant VZ function \big(with respect to $\|\cdot\|\big)$\endslant\ if
$$(f - q) \episum p = 0\ \on\ B.\meqno\VZtwo$$
It follows from (\IOTAfour) and (\VZone) that
$$\hbox{if}\ f\ \hbox{is a VZ function with respect to}\ \|\cdot\|\ \hbox{then}\ \infn_B[f - q] = 0.\meqno\VZthree$$
``VZ'' stands for ``Voisei--Z\u alinescu'', since (\VZtwo) is an extension to Banach SSD spaces of a condition introduced in \cite\VZ, Proposition 3\endcite. 
\defDefinition \Pdef
\medbreak
\noindent
{\bf Definition \Pdef.}\enspace Let $A$ be a subset of a Banach SSD space $B$.   We say that $A$ is \slant $p$--dense\endslant\ if, for all $c \in B$, $\inf p(c - A) = 0$.
 
\medbreak
We now come to our main results on Banach SSD spaces.   Lemma \EXlem(b) is interesting since it tells us that we can determine whether $f$ is a VZ function by inspecting $\P(f)$.
\defLemma \EXlem
\medbreak
\noindent
{\bf Lemma \EXlem.}\enspace\slant Let $B$ be a Banach SSD space and $f \in \PCLSC(B)$.
\par\noindent
{\rm(a)}\enspace  Let $f$ be a VZ function.   Then $\P(f)$ is a $q$--positive subset of $B$ and
$$c \in B \qlr \dist(c,\P(f)) \le \rttwo\sqrt{(f - q)(c)}.\meqno\EXone$$
{\rm(b)}\enspace  The following three conditions are equivalent:
\par
{\rm(i)}\enspace $f$ is a VZ function.
\par
{\rm(ii)}\enspace $f \ge q$ on $B$\quad and, for all $c \in B$ there exists a bounded sequence $\{a_n\}_{n \ge 1}$ of elements of $\P(f)$ such that $\limn_{n \to \infty}p(c - a_n) = 0$.
\par
{\rm(iii)}\enspace $f \ge q$ on $B$ \quand $\P(f)$ is $p$--dense.
\endslant
\Proof (a)\enspace (\VZthree) implies that\quad $f \ge q$ on $B$,\quad and so $\P(f)$ is defined.   Since (\EXone) is trivial if $c \in B \setminus \dom\,f$, we can and will suppose that $c \in \dom\,f$.   Let $\eps \in \,]0,1[\,$.   We first prove that there exists a Cauchy sequence $\{b_n\}_{n \ge 1}$ such that, for all $n \ge 1$,
$$(f - q)(b_n) \le (f - q)(c)/16^n \quand \|c - b_n\| \le (1 + \eps)\rttwo\sqrt{(f - q)(c)}.\meqno\EXtwo$$
Since we can take $b_n = c$ if $(f - q)(c) = 0$, we can and will suppose that
$$\alpha := \sqrt{(f - q)(c)} > 0.\meqno\EXthree$$
Let $\lambda := \eps/(3 + \eps) \in \,]0,1/4[\,$ and write $b_0 := c$.   Then we can choose inductively $b_1,b_2,\dots \in B$ such that, for all $n \ge 1$,\quad $(f - q)(b_n) + p(b_{n - 1} - b_n) \le \lambda^{2n}\alpha^2$.\quad It follows from this and (\VZthree) that, 
$$\all\ n \ge 1,\qquad p(b_{n - 1} - b_n) \le \lambda^{2n}\alpha^2,\meqno\EXfour$$
and, combining with (\IOTAfour),
$$\all\ n \ge 0,\qquad (f - q)(b_n) \le \lambda^{2n}\alpha^2 \le \alpha^2/16^n.\meqno\EXfive$$
Substituting the first inequality of (\EXfive) into Lemma \ROOTlem, for all $n \ge 1$,
$$-q(b_{n - 1} - b_n) \le \Big[\sqrt{(f - q)(b_{n - 1})} + \sqrt{(f - q)(b_n)}\Big]^2 \le (1 + \lambda)^2\lambda^{2n - 2}\alpha^2.$$
Consequently, since $g(b_{n - 1} - b_n) = p(b_{n - 1} - b_n) - q(b_{n - 1} - b_n)$, (\EXfour) gives,
$$\all\ n \ge 1,\qquad g(b_{n - 1} - b_n) \le (1 + \lambda)^2\lambda^{2n - 2}\alpha^2 + \lambda^{2n}\alpha^2 \le  (1 + 2\lambda)^2\lambda^{2n - 2}\alpha^2,$$
and so, for all $n \ge 1$, $\|b_{n - 1} - b_n\| \le  \rttwo(1 + 2\lambda)\lambda^{n - 1}\alpha$.   Adding up this inequality for $n = 1, \dots, m$, we derive that, for all $m \ge 1$, $\|c - b_m\| \le \rttwo(1 + 2\lambda)\alpha/(1 - \lambda)$.   Since $(1 + 2\lambda)/(1 - \lambda) = 1 + \eps$, this and (\EXfive) give (\EXtwo).   Now set $a = \lim_{n}b_n$, so that $\|c - a\| \le (1 + \eps)\rttwo\sqrt{(f - q)(c)}$.   (\EXfive) and the lower semicontinuity of $f - q$ now imply that $(f - q)(a) \le 0$, that is to say, $a \in \P(f)$.   Since $\dom\,f \ne \emptyset$, it follows that $\P(f) \ne \emptyset$ and so, from Lemma \Llem, $\P(f)$ is a $q$--positive subset of $B$.   We also have
$$\dist(c,\P(f)) \le (1 + \eps)\rttwo\sqrt{(f - q)(c)},$$
and so if we now let $\eps \to 0$, we obtain (\EXone).   This completes the proof of (a).
\smallskip
(b)\enspace Suppose first that (i) is satisfied.   (\VZthree) implies that\quad $f \ge q$ on $B$.\quad Let $c \in B$.   We choose inductively $b_1,b_2,\dots \in B$ such that, for all $n \ge 1$,
$$f(b_n) + g(c - b_n) + q(c) - \xbra{c}{b_n} = (f - q)(b_n) + p(c - b_n) < 1/n^2.$$
Consequently, using (\IOTAfour) and (\IOTAtwo), for all $n \ge 1$,
$$(f - q)(b_n) <  1/n^2,\ p(c - b_n) < 1/n^2\meqno\EXseven$$
and
$$f(b_n) + g(c - b_n) + q(c) - \|\iota\|\|c\|\|b_n\| < 1/n^2.\meqno\EXeight$$
Since $f \in \PCLSC(B)$, $f$ dominates a continuous affine function, and so (\EXeight) and the usual coercivity argument imply that $K := \sup_{n \ge 1}\|b_n\| < \infty$.   From (a) and (\EXseven), there exists $a_n \in \P(f)$ such that $\|a_n - b_n\| \le \rttwo/n$.   Now, from (\IOTAfive), for all $n \ge 1$,
$$\eqalign{|p(c - a_n) - p(c - b_n)|
&\le \half(1 + \|\iota\|)\|a_n - b_n\|(2\|c\| + \|a_n\| + \|b_n\|)\cr
&\le \half(1 + \|\iota\|)\big(2\|c\| + \big(K + \rttwo \big) + K\big)\rttwo/n.}$$
Thus $\lim_{n \to \infty}\big[p(c - a_n) - p(c - b_n)\big] = 0$, and (ii) follows by combining this with (\EXseven).   It is trivial that (ii)$\lr$(iii).   Suppose, finally, that (iii) is satisfied.   Then, for all $c \in B$, 
$$\big((f - q) \episum p\big)(c) \le \infn_{a \in \P(f)} \big[(f - q)(a) + p(c - a)\big] = \inf p(c - \P(f)) = 0,$$
from which\quad $(f - q) \episum p \le 0$ on $B$.\quad   On the other hand, since\quad $f - q \ge 0$ on $B$\quad and, from (\IOTAfour),\quad $p \ge 0$ on $B$,\quad we have\quad $(f - q) \episum p \ge 0$ on $B$.\quad   Thus $f$ is a VZ function, giving (i).\qed
\defLemma \Alem
\medbreak
\noindent
{\bf Lemma \Alem.}\enspace\slant Let $A$ be a closed, $p$--dense and $q$--positive subset of a Banach SSD space $B$.
\smallskip\noindent
{\rm(a)}\enspace   For all $c \in B$, $\inf q(c - A) \le 0$ and $\dist(c,A) \le \rttwo\sqrt{{-}\inf q(c - A)}$.
\smallskip\noindent
{\rm(b)}\enspace Let $h \in \PC(B)$, $h \ge q$ on $B$, and $\P(h) \supset A$.   Then $h$ is a VZ function.
\smallskip\noindent
{\rm(c)}\enspace  $A$ is a maximally $q$--positive subset of $B$.\endslant
\Proof(a)\enspace Let $c \in B$.   Then $\inf g(c - A) + \inf q(c - A) \le \inf p(c - A) = 0$.
Thus
$\half\dist(c,A)^2 = \inf g(c - A) \le - \inf q(c - A)$,
from which (a) is an immediate consequence.
\smallskip
(b)\enspace Clearly, $\P(h)$ is also $p$--dense, and it follows as in Lemma \EXlem(b)\big((iii)$\lr$(i)\big) (which does not use any semicontinuity) that $h$ is a VZ function, which gives (b).
\smallskip
(c)\enspace We suppose that $c \in B$ and $\inf q(c - A) \ge 0$, and we must prove that $c \in A$.   From (a), in fact $\inf q(c - A) = 0$ and $\dist(c,A) = 0$.   Since $A$ is closed, $c \in A$.   This completes the proof of (c).\qed
\defTheorem \EXthm
\medbreak
\noindent
{\bf Theorem \EXthm.}\enspace\slant Let $B$ be a Banach SSD space and $f \in \PCLSC(B)$ be a VZ function.   Then: 
\smallskip\noindent
{\rm(a)}\enspace   For all $c \in B$, $\inf q(c - \P(f)) \le 0$ and $\dist(c,\P(f)) \le \rttwo\sqrt{{-}\inf q(c - \P(f))}$.
\smallskip\noindent
{\rm(b)}\enspace Let $h \in \PC(B)$, $h \ge q$ on $B$, and $\P(h) \supset \P(f)$.   Then $h$ is a VZ function.
\smallskip\noindent
{\rm(c)}\enspace  $\P(f)$ is a maximally $q$--positive subset of $B$.
\smallskip\noindent
{\rm(d)}\enspace $f^@ \in \PCLSC(B)$, $f^@$ is a VZ function and $\P\big(f^@\big) = \P(f)$.\endslant
\Proof (a), (b) and (c) are immediate from Lemma \EXlem(b)\big((i)$\lr$(iii)\big) and the corresponding parts of Lemma \Alem.
\smallskip
(d)\enspace Let $c \in B$. Then, since\quad $q \le p$ on $B$,\quad Definition \FATdef\ gives
$$q(c) - f^@(c) = \infn_{b \in B}\big[f(b) - \xbra{b}{c} + q(c)\big] = \big((f - q) \episum q\big)(c) \le \big((f - q) \episum p\big)(c)= 0,$$
and so\quad $f^@ \ge q$ on $B$.\quad   It now follows from Lemma \FFATlem(b) that $\P\big(f^@\big) \supset \P(f)$, and so (b) and (c) imply that $f^@$ is a VZ function and $\P\big(f^@\big) = \P(f)$.   Since $\P(f) \ne \emptyset$, it is evident that $f^@ \in \PCLSC(B)$.   ($\P(f)$ is closed because $f$ is lower semicontinuous.)\qed
\defRemark \MLTrem
\medbreak
\noindent
{\bf Remark \MLTrem.}\enspace In general, Theorem \EXthm(a) is strictly stronger than Lemma \EXlem(a).   While this can be proved directly, we will see in Remark \FPHIrem\ that it follows easily from the properties of the $\Phi$--functions.   We will also see in Remark \FPHIrem\ that the constant $\rttwo$ in (\EXone) is sharp.
\medbreak
The proof of Theorem \EXthm\ relies heavily on the lower semicontinuity of $f$.   We will show in Corollary \Hcor\ below that part of Theorem \EXthm(d) can be recovered even if $f$ is not assumed to be lower semicontinuous.
\defCorollary \Hcor
\medbreak
\noindent
{\bf Corollary \Hcor.}\enspace\slant Let $B$ be a Banach SSD space and $f \in \PC(B)$ be a VZ function.   Then $f^@ \in \PCLSC(B)$, $f^@$ is a VZ function and $\P\big(f^@\big)$ is  a maximally $q$--positive subset of $B$.\endslant
\Proof Let $\fbar$ be the lower semicontinuous envelope of $f$.   Since $q$ is continuous and\quad $f \ge q$ on $B$,\quad it follows that\quad $f \ge \fbar \ge q$ on $B$.\quad   Thus, from  (\IOTAfour),
$$0 = (f - q) \episum p \ge (\fbar - q) \episum p \ge 0 \episum p = 0\ \on\ B,$$      
and so $\fbar$ is a VZ function.   Since $\fbar \in \PCLSC(B)$, Theorem \EXthm(d) implies that $\fbar^@$ is a VZ function also.   It is well known that $\fbar^* = f^*$ on $B^*$ thus, composing with $\iota$ and using (\IOTAsix), $\fbar^@ = f^@$ on $B$. The result now follows from Theorem \EXthm(d,c), with $f$ replaced by $f^@$.\qed
\defDefinition \THdef
\medbreak
\noindent
{\bf Definition \THdef.}\enspace Let $B$ be a Banach SSD space and $A$ be a nonempty $q$--positive subset of $B$.   We define the function $\Theta_A\colon\ B^* \to \rbar$ by:  for all $b^* \in B^*$,
$$\Theta_A(b^*) := \supn_{a \in A}\big[\bra{a}{b^*} - q(a)\big].$$
We define the function $\Phi_A\colon\ B \to \rbar$ by\quad $\Phi_A := \Theta_A \circ \iota$.
\par\noindent
We define the function $^*\Theta_A\colon\ B \to \rbar$ by: for all $c \in B$,
$$^*\Theta_A(c) := \supn_{b^* \in B^*}\big[\bra{c}{b^*} - \Theta_A(b^*)\big].$$
\par
We collect together in Lemma \THlem\ some elementary properties of $\Theta_A$, $\Phi_A$, $^*\Theta_A$, and $\Phi_A^@$.   The properties of $\Phi_A$ and $\Phi_A^@$ have already appeared in \cite\HBM\endcite.
\defLemma \THlem
\medbreak
\noindent
{\bf Lemma \THlem.}\enspace\slant Let $B$ be a Banach SSD space and $A$ be a nonempty $q$--positive subset of $B$.
\smallskip\noindent
{\rm(a)}\enspace  For all $b \in B$,\quad $\Phi_A(b) = \supn_{a \in A}\big[\xbra{a}{b} - q(a)\big] = q(b) - \inf q(b - A)$.
\smallskip\noindent
{\rm(b)}\enspace $\Phi_A \in \PCLSC(B) \quand \Phi_A = q\ \on\ A$.
\smallskip\noindent
{\rm(c)}\enspace $\Theta_A \in \PCLSCST(B^*)$.
\smallskip\noindent
{\rm(d)}\enspace $(^*\Theta_A)^* = \Theta_A$\quad and\quad $(^*\Theta_A)^@ = \Phi_A$.
\smallskip\noindent
{\rm(e)}\enspace $^*\Theta_A \le q$ on $A$.\quad Consequently,\quad $^*\Theta_A \in \PCLSC(B)$.
\smallskip\noindent
{\rm(f)}\enspace $^*\Theta_A \ge {\Phi_A}^@ \ge  \Phi_A \vee q\ \on\ B$. 
\smallskip\noindent
{\rm(g)}\enspace $^*\Theta_A = {\Phi_A}^@ = q\ \on\ A$.
\smallskip\noindent
{\rm(h)}\enspace Let $A$ be maximally $q$--positive.  Then\enspace $^*\Theta_A \ge {\Phi_A}^@ \ge  \Phi_A \ge q$ on $B$\enspace and\enspace $A \subset \P\big(^*\Theta_A\big)$.
\smallskip\noindent
{\rm(i)}\enspace Let $A$ be maximally $q$--positive.  Then\quad $\P\big(^*\Theta_A\big) = \P\big({\Phi_A}^@\big) = \P\big(\Phi_A\big) = A$.\endslant
\Proof(a) is immediate from (\IOTAtwo), (b) from (a), and (c) from (b) and the definition of $\Theta_A$.
\smallskip
The first assertion in (d) follows from (c) and the Fenchel--Moreau theorem for the locally convex space $\big(B^*,w(B^*,B)\big)$, while the second assertion follows from the first by composing with $\iota$, and using (\IOTAsix) and the definition of $\Phi_A$.
\smallskip
(e)\enspace Let $a \in A$.   The definition of $\Theta_A$ implies that, for all $b^* \in B^*$, $\bra{a}{b^*} - \Theta_A(b^*) \le q(a)$.   Taking the supremum over $b^* \in B^*$, $^*\Theta_A(a) \le q(a)$, as required.
\smallskip
(f)\enspace Let $c \in B$.   Then, from (\IOTAtwo), the definition of $\Phi_A$ and (b),
$$\eqalignno{^*\Theta_A(c) 
&\ge \supn_{b \in B}\big[\bra{c}{\iota(b)} - \Theta_A(\iota(b))\big]\cr
&= \supn_{b \in B}\big[\xbra{c}{b} - \Phi_A(b)\big]\quad \big(= {\Phi_A}^@(c)\big)\cr
&\ge \big[\xbra{c}{c} - \Phi_A(c)\big] \vee \supn_{a \in A}\big[\xbra{c}{a} - \Phi_A(a)\big]\cr
&= \big[2q(c) - \Phi_A(c)\big] \vee \supn_{a \in A}\big[\xbra{c}{a} - q(a)\big]\cr
&= \big[2q(c) - \Phi_A(c)\big] \vee \Phi_A(c).}$$
Now if $\Phi_A(c) = \infty$ then obviously $\big[2q(c) - \Phi_A(c)\big] \vee \Phi_A(c) \ge q(c)$, while if $\Phi_A(c) \in \r$ then $\big[2q(c) - \Phi_A(c)\big] \vee \Phi_A(c) \ge \half[2q(c) - \Phi_A(c)\big] + \half\Phi_A(c) = q(c)$.   Thus ${\Phi_A}^@(c) \ge  \Phi(c) \vee q(c)$.   This completes the proof of (f).
\smallskip
(g) is immediate from (e) and (f).
\smallskip
(h)\enspace In this case, for all $b \in B \setminus A$, there exists $a \in A$ such that $q(b - a) < 0$, and so $\inf q(b - A) < 0$.  Thus, from (a),\quad $\Phi_A > q$ on $B \setminus A$.\quad   Combining this with (b),\quad $\Phi_A \ge q$ on $B$\quad and\quad $\P\big(\Phi_A\big) = A$.\quad  Thus (h) follows from (f) and (g).
\smallskip
It is clear from (h) that\quad $A \subset\P\big(^*\Theta_A\big) \subset \P\big({\Phi_A}^@\big) \subset \P\big(\Phi_A\big)$,\quad and so (i) follows from the maximality of $A$.\qed
\defRemark \THrem
\medbreak
\noindent
{\bf Remark \THrem.}\enspace We will see in (\HHATtwo) and (\HHATthree) that $^*\Theta$ and ${\Phi_A}^@$ are both ``upper limiting'' functions in various situations, so the question arises whether these two functions are identical.   If $^*\Theta_A = {\Phi_A}^@$ then $^*\Theta_A$ is obviously $w(B,B)$--lower semicontinuous.   If, conversely, $^*\Theta_A$ is $w(B,B)$--lower semicontinuous then the Fenchel--Moreau\break theorem for the (possibly nonhausdorff) locally convex space $\big(B,w(B,B)\big)$ and Lemma \THlem(d) imply that $^*\Theta_A = \big(^*\Theta_A\big)^{@@} = {\Phi_A}^@$.   The author is grateful to Constantin Z\u alinescu for the following example showing that, in general, the functions $^*\Theta$ and ${\Phi_A}^@$ are not identical.   Let $B$ be a Banach space,\quad $\xbra\cdot\cdot = 0$ on $B \times B$\quad and $A$ be a nonempty proper closed convex subset of $B$.   Then $^*\Theta_A$ is the indicator function of $A$ and\quad ${\Phi_A}^@ = 0$\quad on $B$.   We do not know what the situation is if $A$ is maximally $q$--positive, or in the special situation of Example \EESPECex.   For the convenience of the reader, we will give a proof of the Fenchel--Moreau theorem for  nonhausdorff locally convex spaces in Theorem \FMthm.
\defTheorem \THthm
\medbreak
\noindent
{\bf Theorem \THthm.}\enspace\slant Let $B$ be a Banach SSD space.
\par\noindent
{\rm(a)}\enspace Let $f \in \PCLSC(B)$,\quad $f \ge q$ on $B$\quad and $A := \P(f) \ne \emptyset$.   Then\quad $^*\Theta_A \ge f \ge \Phi_A$ on $B$\quad and\quad ${\Phi_A}^* \ge f^* \ge \Theta_A$ on $B^*$.
\smallskip\noindent
{\rm(b)}\enspace Let $A$ be a maximally $q$--positive subset of $B$, $h \in \PC(B)$ and\quad $^*\Theta_A \ge h \ge \Phi_A$ on $B$.\quad   Then\quad $h \ge q$ on $B$,\quad $h^@ \ge q$ on $B$\quad and\quad $\P(h) = \P\big(h^@\big) = A$.
\smallskip\noindent
{\rm(c)}\enspace Let $f \in \PCLSC(B)$ be a VZ function and $A := \P(f)$.   Then
$$^*\Theta_A \ge f \ge \Phi_A \ge q\ \on\ B \quand {\Phi_A}^* \ge f^* \ge \Theta_A\ \on\ B^*.\meqno\THone$$
Now let $h \in \PC(B)$ \quand $^*\Theta_A \ge h \ge \Phi_A$ on $B$.   Then $h$ and $h^@$ are VZ functions.   In particular, $\P\big(^*\Theta_A\big) = \P\big({\Phi_A}^@\big) = \P\big(\Phi_A\big) = \P(f)$ and $\Phi_A$, ${\Phi_A}^@$ and $^*\Theta_A$ are all VZ functions.
\endslant
\Proof(a)\enspace Let $b \in B$ and $a \in \P(f)$.   Then, from Lemma \FFATlem(a), $f(b)\ge \xbra{b}{a} - q(a)$.   Taking the supremum over $a \in \P(f)$ and using Lemma \THlem(a), $f(b) \ge \Phi_A(b)$.   Thus\quad $f \ge \Phi_A$ on $B$\quad and, taking conjugates,\quad ${\Phi_A}^* \ge f^*$ on $B^*$.\quad   Now, for all $b^* \in B^*$,
$$\eqalign{f^*(b^*)
&= \supn_{b \in B}\big[\bra{b}{b^*} - f(b)\big]
\ge \supn_{a \in \P(f)}\big[\bra{a}{b^*} - f(a)\big]\cr
&= \supn_{a \in \P(f)}\big[\bra{a}{b^*} - q(a)\big] = \Theta_A(b^*).}$$
Thus\quad $f^* \ge \Theta_A$ on $B^*$.\quad   Taking conjugates and using the Fenchel--Moreau theorem for the normed space $B$,\quad $^*\Theta_A \ge f$ on $B$.\quad   This completes the proof of (a).
\smallskip
(b)\enspace From Lemma \THlem(h),
$$^*\Theta_A \ge h \ge \Phi_A \ge q\ \on\ B,\quad \hbox{from which}\quad \P\big(^*\Theta_A\big) \subset \P(h) \subset \P\big(\Phi_A\big).\meqno\HHATtwo$$
It is clear from our assumptions that ${\Phi_A}^@ \ge h^@ \ge (^*\Theta_A)^@$ on $B$.   If we now combine this with Lemma \THlem(d,h), we derive that
$${\Phi_A}^@ \ge h^@ \ge \Phi_A \ge q\ \on\ B,\quad  \hbox{from which}\quad \P\big({\Phi_A}^@\big) \subset \P\big(h^@\big) \subset \P\big(\Phi_A\big).\meqno\HHATthree$$
(b) now follows from (\HHATtwo), (\HHATthree) and Lemma \THlem(i).
\smallskip
(c) The assertions about $f$ follow from (\VZthree), Theorem \EXthm(c), (a) and Lemma \THlem(h), the assertions about $h$ and $h^@$ follow from Theorem \EXthm(c,b) and (b), and then the assertions about $\Phi_A$, ${\Phi_A}^@$ and $^*\Theta_A$ follow from Theorem \EXthm(c) and Lemma \THlem(h,i).\qed
\medbreak
In Theorem \MAXthm\ below, we show that $^*\Theta_A$ has a certain maximal property.   This result was motivated by results originally proved by Burachik and Svaiter in \cite\BS\endcite\ for maximally monotone multifunctions.
\defTheorem \MAXthm
\medbreak
\noindent
{\bf Theorem \MAXthm.}\enspace\slant Let $A$ be a nonempty $q$--positive subset of a Banach SSD space $B$ and
$$\sigma_A := \sup\big\{h\colon\ h \in \PCLSC(B),\ h \le q\ \on\ A\big\}.$$
Then\quad $^*\Theta_A = \sigma_A$ on $B$.\endslant
\Proof Let $h \in \PCLSC(B)$ and $h \le q$ on $A$.   The Fenchel Young inequality implies that, for all $b^* \in B^*$ and $a \in A$, $h^*(b^*) \ge \bra{a}{b^*} - h(a) \ge \bra{a}{b^*} - q(a)$.   Thus, taking the supremum over $a \in A$, $h^*(b^*) \ge \Theta_A(b^*)$.   In other words,\quad $h^* \ge \Theta_A$ on $B^*$.\quad Taking conjugates and using the Fenchel Moreau theorem for the normed space $B$,\quad $^*\Theta_A \ge h$ on $B$.\quad   It follows by taking the supremum over $h$ that\quad $^*\Theta_A \ge \sigma_A$ on $B$.\quad   On the other hand, it is clear from Lemma \THlem(e) that\quad $\sigma_A \ge {^*\Theta_A}\ \on\ B$.
\qed
\defRemark \FPHIrem
\medbreak
\noindent
{\bf Remark \FPHIrem.}\enspace Let $B$ be a Banach SSD space and $f \in \PCLSC(B)$ be a VZ function.   We know from  Theorem \THthm(c) that $\P\big(\Phi_{\P(f)}\big) = \P(f)$, $\Phi_{\P(f)}$ is a VZ function and\quad $\Phi_{\P(f)} \le f$ on $B$.\quad Thus Lemma \EXlem(a) implies that, for all $c \in B$,
$$\dist(c,\P(f)) = \dist\big(c,\P\big(\Phi_{\P(f)}\big)\big) \le \rttwo\sqrt{\big(\Phi_{\P(f)} - q\big)(c)} \le \rttwo\sqrt{(f - q)(c)}.$$  
From Lemma \THlem(a), $\big(\Phi_{\P(f)} - q\big)(c) = \Phi_{\P(f)}(c) - q(c) = -\inf q(c - \P(f))$, thus we have
$$\dist(c,\P(f)) \le \rttwo\sqrt{-\inf q(c - \P(f))} \le \rttwo\sqrt{(f - q)(c)}.$$
This shows that Theorem \EXthm(a) is as least as strong as Lemma \EXlem(a).   Now let $E := \r$ and $B$ be the Banach SSD space $\r^2$ as in Example \EESPECex, using the norm $\|\cdot\|_{2,1}$.   Define $f \in \PCLSC(B)$ by $f(x_1,x_2) := \half(x_1^2 + x_2^2)$.   Then $(f - q)(x_1,x_2) = \half(x_1^2 + x_2^2) - x_1x_2 = \half(x_1 - x_2)^2$
and
$p(x_1,x_2) = \half(x_1^2 + x_2^2) + x_1x_2 = \half(x_1 + x_2)^2$.
Let $c := (z_1,z_2) \in B$ and $b := \big(\half(z_1 + z_2),\half(z_1 + z_2)\big) \in B$.   Then $(f - q)(b) = 0$ and $p(c - b) = 0$.   Consequently, $f$ is a VZ function.  Now $\P(f)$ is the diagonal of $\r^2$ and so, by direct computation, for all $c = (x_1,x_2) \in \r^2$,
$-\inf q(c - \P(f)) = \fourth(x_1 - x_2)^2$.   Since $\fourth(x_1 - x_2)^2 < \half(x_1 - x_2)^2$ when $x_1 \ne x_2$, Theorem \EXthm(a) is strictly stronger than Lemma \EXlem(a) in this case.\par
Now let $h := \Phi_{\P(f)}$.   Lemma \THlem(a) gives us that, for all $(x_1,x_2) \in B$,
$$\sqrt{(h - q)(x_1,x_2)} = \sqrt{\fourth(x_1 - x_2)^2} = \half|x_1 - x_2|.$$
On the other hand, by direct computation, $\dist\big((x_1,x_2),\P(h)\big) = \rthalf|x_1 - x_2|$.   Thus the constant $\rttwo$ in (\EXone) is sharp.   The genesis of this argument and example can be found in the results of Mart\'\i nez-Legaz and Th\'era in \cite\MT\endcite.
\defRemark \NOQTrem
\medbreak
\noindent
{\bf Remark \NOQTrem.}\enspace We note that the inequalities for $B$ in (\THone) have four functions, while the inequality for $B^*$ has only three.   The reason for this is that we do not have a function on $B^*$ that plays the role that the function $q$ plays on $B$.   We will introduce such a function in Definition \DUALdef.
\defSection \BSSDAPPSsec
\bigbreak   
\centerline{\bf \BSSDAPPSsec\quad Applications of Section \BSSDsec\ to $E \times E^*$}
\medskip\noindent
In this section, we suppose that $E$ is a nonzero Banach space, and follow the notation of Example \EENex.   Let $A$ be a nonempty monotone subset of $E \times E^*$.   In this case, the definitions and results obtained in Definition \THdef\ and Lemma \THlem\ specialize as follows. The function $\Theta_A \in \PCLSCST(E\dbs \times E^*)$ is defined by:
$$\Theta_A(x\dbs,x^*) := \supn_{(s,s^*) \in A}\big[\bra{s}{x^*} + \bra{s^*}{x\dbs} - \bra{s}{s^*}\big].$$
The function $\Phi_A \in \PCLSC(E \times E^*)$ is defined by:
$$\Phi_A(x,x^*) = \supn_{(s,s^*) \in A}\big[\bra{x}{s^*} + \bra{s}{x^*} - \bra{s}{s^*}\big].$$
$\Phi_A$ is the \slant Fitzpatrick function\endslant\ of $A$, first introduced in \cite\FITZ\endcite, which has been discussed by many authors in recent years.   The function $^*\Theta_A \in \PCLSC(E \times E^*)$ is defined by:
$$^*\Theta_A(y,y^*) := \supn_{(x\dbs,x^*) \in E\dbs \times E^*}\big[\bra{y}{x^*} + \bra{y^*}{x\dbs} - \Theta_A(x\dbs,x^*)\big].$$
Then $(^*\Theta_A)^* = \Theta_A$ and $(^*\Theta_A)^@ = \Phi_A$.   Furthermore,
$$^*\Theta_A \ge {\Phi_A}^@ \ge  \Phi_A \vee q\ \on\ E \times E^* \quand ^*\Theta_A = {\Phi_A}^@ = \Phi_A = q\ \on\ A.$$
If $f \in \PC(E \times E^*)$ and\quad $f \ge q$ on $ E \times E^*$\quad   then we define $\M f$ to be the monotone set $\{(x,x^*) \in E \times E^*\colon f(x,x^*) = \bra{x}{x^*}\}$.   $\M f$ is identical with $\P(f)$ as in Definition \POSFdef, but the ``$\M$'' notation seems more appropriate in this case.   Continuing with the consequences of Lemma \THlem, we have:
$$^*\Theta_A \ge {\Phi_A}^@ \ge  \Phi_A \ge q\ \on\ E \times E^*\quand \M\big(^*\Theta_A\big) = \M\big({\Phi_A}^@\big) = \M\big(\Phi_A\big) = A.$$
\par
\medbreak
The following results are then immediate from Theorems \THthm\ and \MAXthm.   The expression $\sup\big\{h\colon\ h \in \PCLSC(E \times E^*),\ h \le q\ \on\ A\big\}$ that appears in Theorem \NTHthm(b) was first introduced by Burachik and Svaiter in \cite\BS\endcite\ (for $A$ maximally monotone) and further studied by Marques Alves and Svaiter in \cite\ASTHREE\endcite.   The analysis of Lemma \THlem\ and Theorem \THthm\ suggests that the natural framework in which to consider these results is that of Banach SSD spaces.
\defTheorem \NTHthm
\medbreak
\noindent
{\bf Theorem \NTHthm.}\enspace\slant Let $E$ be a nonzero Banach space, $E \times E^*$ be normed as in Example \EENex, and $A$ be a nonempty monotone subset of $E \times E^*$.
\par\noindent
{\rm(a)}\enspace Let $f \in \PCLSC(E \times E^*)$, $f \ge q$ on $E \times E^*$ and $A := \M f \ne \emptyset$.   Then
$$^*\Theta_A \ge f \ge \Phi_A\ \on\ E \times E^*,\quand {\Phi_A}^* \ge f^* \ge \Theta_A\ \on\ E\dbs \times E^*.$$
{\rm(b)}\enspace Let $A$ be maximally monotone, $h \in \PC(E \times E^*)$ and\quad $^*\Theta_A \ge h \ge \Phi_A$ on $E \times E^*$.\quad Then\quad $h \ge q$ on $E \times E^*$,\quad $h^@ \ge q$ on $E \times E^*$ \quand $\M h = \M\big(h^@\big) = A$.
\smallskip\noindent
{\rm(c)}\enspace Let $f \in \PCLSC(E \times E^*)$ be a VZ function and $A := \M f$.   Then
$$^*\Theta_A \ge f \ge \Phi_A \ge q\ \on\ E \times E^* \quand {\Phi_A}^* \ge f^* \ge \Theta_A\ \on\ E\dbs \times E^*$$
Now let $h \in \PC(E \times E^*)$ \quand $^*\Theta_A \ge h \ge \Phi_A$ on $E \times E^*$.\quad Then $h$ and $h^@$ are VZ functions.   In particular, $\M\big(^*\Theta_A\big) = \M\big({\Phi_A}^@\big) = \M\big(\Phi_A\big) = A$, and $\Phi_A$, ${\Phi_A}^@$ and $^*\Theta_A$ are all VZ functions.
$$^*\Theta_A = \sup\big\{h\colon\ h \in \PCLSC(E \times E^*),\ h \le q\ \on\ A\big\}.\leqno{\rm(d)}$$
\endslant
\defSection \DUALsec
\medbreak   
\centerline{\bf \DUALsec\quad Banach SSD dual spaces}
\defDefinition \DUALdef
\medbreak
\noindent
{\bf Definition \DUALdef.}\enspace Let $(B,\|\cdot\|)$ be a Banach SSD space and $(B^*,\|\cdot\|)$ be the norm--dual of $B$.   We say that $(B^*,\ybra{\cdot}{\cdot})$ is a \slant Banach SSD dual of $B$\endslant\ if $\ybra{\cdot}{\cdot}\colon\ B^* \times B^* \to \r$ is a symmetric bilinear form,
$$\all\ b \in B\ \and\ c^* \in B^*,\qquad \ybra{\iota(b)}{c^*} = \bra{b}{c^*}.\meqno\DUALone$$
Writing $\qt(c^*) := \half\ybra{c^*}{c^*}$ and $\pt(c^*) := \half\|c^*\|^2 +\qt(c^*)$, we suppose also that
$$\pt \ge 0\ \on\ B^*.\meqno\DUALtwo$$
Now if we take $c^* = \iota(c)$ in (\DUALone) and use (\IOTAtwo), we obtain
$$\all\ b,c\in B,\qquad \ybra{\iota(b)}{\iota(c)} = \Bra{b}{\iota(c)} = \xbra{b}{c},\meqno\DUALthree$$
from which
$$\qt \circ \iota = q.\meqno\DUALfour$$
It is easy to see from these definitions that,
$$\all\ b^* \in B^*,\qquad \Theta_A(b^*) = \qt(b^*) - \inf\qt(b^* - \iota(A)).\meqno\DUALfive$$
This should be compared with Lemma \THlem(a).
\defDefinition \PSEUDOdef
\medbreak
\noindent
{\bf Definition \PSEUDOdef.}\enspace Let $(B,\|\cdot\|)$ be a Banach SSD space and $(B^*,\ybra{\cdot}{\cdot})$ be a Banach SSD dual of $B$.
We say that $\iota(B)$ is \slant $\pt$--dense in\endslant\ $B^*$ if
$$\all\ b^* \in B^*,\quad \inf\pt\big(b^* - \iota(B)\big) = 0.\meqno\PSEUDOone$$\par
\defRemark \DUALrem
\medbreak
\noindent
{\bf Remark \DUALrem.}\enspace In Example \Hex\ with $\|T\| \le 1$ (see also Remark \NORMrem), for all $c \in B$, $\iota(c) = Tc$.   Suppose now that $T^2$ is the identity on $B$.   Since $B^* = B$,
$$\all\ b \in B\ \and\ c^* \in B^* = B,\qquad \xbra{\iota(b)}{c^*} =  \xbra{Tb}{c^*} = \bra{T^2b}{c^*} = \bra{b}{c^*}.$$
Thus (\DUALone) is satisfied with $\ybra{\cdot}{\cdot} := \xbra{\cdot}{\cdot}$, and so \big($B,\xbra{\cdot}{\cdot}$\big) is its own Banach SSD dual.   We note that $T^2$ is the identity on $B$ in (a), (b) and (c) of Example \Hex.
\defExample \EEDUALex
\medbreak
\noindent
{\bf Example \EEDUALex.}\enspace We now continue our discussion of Examples \EEex, \EENex\ and \EESPECex.   We recall that $B = E \times E^*$,  $B^* = E\dbs \times E^*$ and, for all $(x,x^*) \in B$, $\iota(x,x^*) = (\wh{x},x^*)$.   We define the symmetric bilinear form $\ybra{\cdot}{\cdot}\colon\ B^* \times B^* \to \r$ by
$$\ybra{b^*}{c^*} := \bra{y^*}{x\dbs} + \bra{x^*}{y\dbs}\quad\big(b^* = (x\dbs,x^*) \in B^*,\ c^* = (y\dbs,y^*) \in B^*\big).$$
It is then easily checked from (\EEone) that (\DUALone) is satisfied and, for all $c^* = (y\dbs,y^*) \in B^*$, $\qt(c^*) = \half\big[\bra{y^*}{y\dbs} + \bra{y^*}{y\dbs}\big] = \bra{y^*}{y\dbs}$.   We now discuss briefly the limitations of this definition.   Let $E := \r$ and $B$ be the SSD space $\r^2$ as in Example \EEex, using the norm $2\|\cdot\|_{2,1}$.   As we observed in Example \EESPECex, $\r^2$ is a Banach SSD space under $\|\cdot\|_{2,1}$, and consequently also a Banach SSD space under the larger norm $2\|\cdot\|_{2,1}$.  Since $\iota$ is the identity on $\r^2$, (\DUALthree) implies that $\ybra{\cdot}{\cdot} := \xbra{\cdot}{\cdot}$.   Now the norm on $B^* = B$ dual to $2\|\cdot\|_{2,1}$ is $\half\|\cdot\|_{2,1}$.  Since $\pt(1,-1) = \f18(2) + (1)(-1) = -\f34 < 0$, $B$ does not admit a Banach SSD dual.   We now return to the general case.   If $c^* = (y\dbs,y^*) \in B^*$ then
$$\half\|c^*\|_{1,\tau}^2 + \qt(c^*) \ge \fourth\big(\tau\|y\dbs\| + \|y^*\|/\tau\big)^2 - \|y\dbs\|\|y^*\| = \fourth\big(\tau\|y\dbs\| - \|y^*\|/\tau\big)^2 \ge 0.$$
Consequently, $\big(B^*,\|\cdot\|_{1,\tau}\big)$ is a Banach SSD dual of $\big(B,\|\cdot\|_{\infty,\tau}\big)$.   Since\quad $\|\cdot\|_{\infty,\tau} \ge \|\cdot\|_{2,\tau} \ge \|\cdot\|_{1,\tau}$ on $B^*$,\quad $\big(B^*,\|\cdot\|_{2,\tau}\big)$ is a Banach SSD dual of $\big(B,\|\cdot\|_{2,\tau}\big)$ and $\big(B^*,\|\cdot\|_{\infty,\tau}\big)$ is a Banach SSD dual of $\big(B,\|\cdot\|_{1,\tau}\big)$.   Next, if $b^* = (y\dbs,y^*) \in B^*$ and $\eps > 0$ then there exists $z^* \in E^*$ such that $\|z^*\| \le \|\tau y\dbs\|$ and $\bra{z^*}{\tau y\dbs} \ge \|\tau y\dbs\|^2 - \eps$.   Let $c := (0,y^* + \tau z^*) \in B$, so that $b^* - \iota(c) = (y\dbs,-\tau z^*) \in B^*$.   Thus
$$\eqalign{\half\|b^* - \iota(c)\|_{\infty,\tau}^2 &+ \qt(b^* - \iota(c)) = \big(\tau\|y\dbs\| \vee \|z^*\|\big)^2 - \bra{\tau z^*}{y\dbs}\cr
&= \big(\|\tau y\dbs\| \vee \|z^*\|\big)^2 - \bra{z^*}{\tau y\dbs} = \|\tau y\dbs\|^2 - \bra{z^*}{\tau y\dbs} \le \eps.}$$
Consequently, if $B$ is normed by $\|\cdot\|_{1,\tau}$ then $\iota(B)$ is $\pt$--dense in $B^*$.   Since\quad $\|\cdot\|_{1,\tau} \le \|\cdot\|_{2,\tau} \le \|\cdot\|_{\infty,\tau}$ on $B^*$,\quad the same is true if $B$ is normed by $\|\cdot\|_{2,\tau}$ or $\|\cdot\|_{\infty,\tau}$.
\medskip
\medbreak
We now recall Rockafellar's formula for the conjugate of a sum:
\defLemma \RSUMlem
\medbreak
\noindent
{\bf Lemma \RSUMlem.}\enspace\slant Let $X$ be a nonzero real Banach space and $f \in \PC(X)$, and let $h \in \PC(X)$ be real--valued and continuous.   Then, for all $x^* \in X^*$,
$$(f + h)^*(x^*) = \minn_{y^* \in X^*}\big[f^*(y^*) + h^*(x^* - y^*)\big].$$\par
\endslant
\Proof  See Rockafellar, \cite\FENCHEL, Theorem 3(a), p.\ 85\endcite, Z\u alinescu, \cite\ZBOOK, Theorem 2.8.7(iii), p.\ 127\endcite, or \cite\HBM, Corollary 10.3, p.\ 52\endcite.\qed
\defRemark \SHARPrem
\medbreak
\noindent
{\bf Remark \SHARPrem.}\enspace \cite\HBM, Theorem 7.4, p.\ 43\endcite\ contains a version of the Fenchel duality theorem with a sharp lower bound on the functional obtained.
\medbreak
Our next result exhibits a certain pleasing symmetry between $B$ and $B^*$.  
\defLemma \SSDDlem
\medbreak
\noindent
{\bf Lemma \SSDDlem.}\enspace\slant Let $B$ be a Banach SSD space with a Banach SSD dual $B^*$ and $f \in \PC(B)$.   Then\quad $\big((f - q) \episum p\big) + \big((f^* - \qt) \episum \pt\big)\circ\iota = 0$  on  $B$.\endslant
\Proof Let $c \in B$.   Define $h\colon\ B \to \r$ by $h(b) := g(c - b)$.   Then, by direct computation using the fact that $g$ is an even function,
$$\all\ c^* \in B^*,\qquad h^*(c^*) = g^*(c^*) + \bra{c}{c^*}.\meqno\Hone$$
Then, using (\IOTAtwo), the continuity of $h$, Lemma \RSUMlem, (\Hone), (\DUALfour) and the fact that, for all $c^* \in B^*$, $g^*(c^*) = \half\|c^*\|^2$,
$$\eqalign{-\big((f - q) \episum p\big)(c)
&= \supn_{b \in B}\big[-(f - q)(b) - p(c - b)\big]\cr
&= \supn_{b \in B}\big[\bra{b}{\iota(c)} -f(b) - h(b)\big] - q(c) = (f + h)^*\big(\iota(c)\big) - q(c)\cr
&=\minn_{b^* \in B^*}\big[f^*(b^*) + h^* \big(\iota(c) - b^*\big)\big] - q(c)\cr
&= \minn_{b^* \in B^*}\big[f^*(b^*) + g^*\big(\iota(c) - b^*\big) + \Bra{c}{\iota(c) - b^*}\big] - q(c)\cr
&= \minn_{b^* \in B^*}\big[f^*(b^*) + g^*\big(\iota(c) - b^*\big) - \ybra{\iota(c)}{b^*} + \qt\big(\iota(c)\big)\big]\cr
&= \minn_{b^* \in B^*}\big[(f^* - \qt)(b^*) + \pt\big(\iota(c) - b^*\big)\big]\cr
&= \big((f^* - \qt) \episum \pt\big)\big(\iota(c)\big).}$$
This completes the proof of Lemma \SSDDlem.\qed    
\defDefinition \MASdef
\medbreak
\noindent
{\bf Definition \MASdef.}\enspace Let $B$ be a Banach SSD space with Banach SSD dual $B^*$ and $f \in \PC(B)$.   We say that $f$ is an \slant MAS function\endslant\ if \quad $f \ge q$ on $B$\quad and \quad $f^* \ge \qt$ on $B^*$.\quad   This is an extension to Banach SSD spaces of the concept introduced by Marques Alves and Svaiter in \cite\ASTWO\endcite\ for the situation described in Example \EEDUALex.
\defTheorem \EQthm
\medbreak
\noindent
{\bf Theorem \EQthm.}\enspace\slant Let $B$ be a Banach SSD space with Banach SSD dual $B^*$ and $f \in \PC(B)$.
\par\noindent
{\rm(a)}\enspace Let $f$ be an MAS function.   Then $f$ is a VZ function.
\par\noindent
{\rm(b)}\enspace Let $\iota(B)$ be $\pt$--dense in $B^*$ and $f$ be a VZ function.   Then $f$ is an MAS function.\par\noindent
{\rm(c)}\enspace Let $\iota(B)$ be $\pt$--dense in $B^*$.   Then $f$ is a VZ function if, and only if, $f$ is an MAS function.\endslant
\Proof(a)\enspace We have \big(using (\IOTAfour)\big)\quad $f - q \ge 0$ and $p \ge 0$ on $B$,\quad and \big(using (\DUALtwo)\big),\quad $f^* - \qt \ge 0$ and $\pt \ge 0$ on $B^*$.\quad   Thus\quad $(f - q) \episum p \ge 0$ and $\big((f^* - \qt) \episum \pt\big)\circ\iota \ge 0$  on $B$.\quad  It now follows from Lemma \SSDDlem\ that $f$ is a VZ function.\par
(b)\enspace Let $b^* \in B^*$ and $c \in B$.   Then, from Lemma \SSDDlem\ again,
$$(f^* - \qt)(b^*) + \pt\big(\iota(c) - b^*\big)
\ge \big((f^* - \qt) \episum \pt\big)\big(\iota(c)\big)
= -\big((f - q) \episum p\big)(c) = 0.$$
Taking the infimum over $c \in B$ and using (\PSEUDOone),\quad $(f^* - \qt)(b^*) \ge 0$ on $B^*$.\quad   Since this holds for all $b^* \in B^*$, $f$ is an MAS function.
\par
(c) is immediate from (a) and (b).\qed
\medbreak
In Theorem \MULTIthm, we shift the emphasis from the properties of a given function $f \in \PC(B)$ to the properties of a given maximally $q$--positive subset $A$ of $B$.   We note that (a), (b), (c), (f) and (g) of Theorem \MULTIthm\ do not involve any functions on $B$ other than those introduced in Definition \THdef.
\defTheorem \MULTIthm
\medbreak
\noindent
{\bf Theorem \MULTIthm.}\enspace\slant Let $(B,\|\cdot\|)$ be a Banach SSD space with Banach SSD dual $B^*$ and $\iota(B)$ be $\pt$--dense in $B^*$.   Let $A$ be a maximally $q$--positive subset of $B$.   Then the following conditions are equivalent:
\par\noindent
{\rm(a)}\enspace For all $b^* \in B^*$,\quad $\inf \qt\big(b^* - \iota(A)\big) \le 0$.
\smallskip\noindent
{\rm(b)}\enspace $\Theta_A \ge \qt$ on $B^*$.
\smallskip\noindent
{\rm(c)}\enspace ${\Phi_A}^* \ge \qt$ on $B^*$.
\smallskip\noindent
{\rm(d)}\enspace There exists an MAS function $f \in \PCLSC(B)$ such that $\P(f) = A$.
\smallskip\noindent
{\rm(e)}\enspace There exists a VZ function $f \in \PCLSC(B)$ such that $\P(f) = A$.
\smallskip\noindent
{\rm(f)}\enspace $\Phi_A$ is a VZ function.
\smallskip\noindent
{\rm(g)}\enspace $^*\Theta_A$ is a VZ function.
\smallskip\noindent
{\rm(b$_1$)}\enspace If $h \in \PC(B)$ and\quad $^*\Theta_A \ge h$ on $B$\quad then\quad $h^* \ge \qt$ on $B^*$.
\smallskip\noindent
{\rm(b$_2$)}\enspace If $h \in \PCLSC(B)$ and\quad $^*\Theta_A \ge h \ge \Phi_A$ on $B$\quad then\quad $h^* \ge \qt$ on $B^*$.
\smallskip\noindent
{\rm(c$_1$)}\enspace There exists $h \in \PCLSC(B)$ such that\enspace $^*\Theta_A \ge h \ge \Phi_A$ on $B$\enspace and\enspace $h^* \ge \qt$ on $B^*$.
\smallskip\noindent
{\rm(c$_2$)}\enspace There exists $h \in \PC(B)$ such that\enspace $h \ge \Phi_A$ on $B$\enspace and\enspace $h^* \ge \qt$ on $B^*$.
\endslant
\Proof The equivalence of (a) and (b) is immediate from (\DUALfive). Taking the conjugate of the inequality in Lemma \THlem(f) and using Lemma \THlem(d) implies that\quad ${\Phi_A}^* \ge \Theta_A$ on $B^*$.\quad   Thus (b)$\lr$(c).   If (c) is satisfied then Lemma \THlem(b,h,i) give (d) with $f := \Phi_A$.   It is immediate from Theorem \EQthm(c) that (d)$\lr$(e).   If (e) is satisfied then Theorem \THthm(c) gives (f) and (g).   If (f) or (g) is satisfied then, from Theorem \EQthm(c) again, $\Phi_A$ or $^*\Theta_A$ (respectively) are MAS functions.   The first of these possibilities implies (c), and the second of these possibilities together with Lemma \THlem(d) implies (b).   Thus (a), (b), (c), (d), (e), (f) and (g) are equivalent.
\par
If $h \in \PC(B)$ and\quad  $^*\Theta_A \ge h$ on $B$\quad  then, from Lemma \THlem(d),\quad  $h^* \ge (^*\Theta_A)^* = \Theta_A$\  on $B^*$,\quad thus (b) implies (b$_1$).   It is trivial that (b$_1$) implies (b$_2$), and it follows by taking $h := {^*\Theta_A}$ and using Lemma \THlem(f,d) that (b) is true.   Thus (b), (b$_1$) and (b$_2$) are equivalent.
\par
If (c) is true then (c$_1$) follows by taking $h := \Phi_A$ and using Lemma \THlem(b,f).   It is trivial that (c$_1$) implies (c$_2$).   If $h \in \PC(B)$ and \quad $h \ge \Phi_A$ on $B$\quad then\quad ${\Phi_A}^* \ge h^*$ on $B^*$, and (c$_2$) implies (c).   Thus (c), (c$_1$) and (c$_2$) are equivalent.\qed 
\defSection \DUALAPPSsec
\bigbreak   
\centerline{\bf \DUALAPPSsec\quad Applications of Section \DUALsec\ to $E \times E^*$}
\medskip\noindent
In this section, we suppose that $E$ is a nonzero Banach space, and show how the results of Section \DUALsec\ can be applied to Example \EEDUALex.   We refer the reader to Section \BSSDAPPSsec\ for the definitions of $\Theta_A$ and $\Phi_A$ in this case.
\defRemark \DIFFrem
\medbreak
\noindent
{\bf Remark \DIFFrem.}\enspace Before proceeding with our analysis, we make some remarks about the essential difference between the concepts of MAS function introduced in Definition \MASdef\ and VZ function introduced in Definition \VZdef.   As observed in Example \EEDUALex, we have $(E \times E^*)^* = E\dbs \times E^*$ and $\qt(x\dbs,x^*) = \bra{x^*}{x\dbs}$, so we have all the information needed to decide whether a function $f \in \PC(E \times E^*)$ is an MAS function.   The situation with VZ functions is different since that involves the function $g$ in an essential way, and this is determined by the precise norm we are using on $E \times E^*$.   In order to clarify the situation, we make the following definition.
\defDefinition \SPECdef
\medbreak
\noindent
{\bf Definition \SPECdef.}\enspace We say that the norm $\|\cdot\|$ on $E \times E^*$ is \slant special\endslant\ if, for some $\tau > 0$, $\|\cdot\|$ is identical with one of the norms $\|\cdot\|_{1,\tau}$, $\|\cdot\|_{2,\tau}$ or $\|\cdot\|_{\infty,\tau}$ introduced in Example \EESPECex.   As we pointed out in the comments in Example \EEDUALex, if $E \times E^*$ is normed by a special norm then $\iota(E \times E^*)$ is $\pt$--dense in $(E \times E^*)^*$.    
\defTheorem \INVARthm
\medbreak
\noindent
{\bf Theorem \INVARthm.}\enspace\slant Let $E$ be a nonzero Banach space and $f \in \PC(E \times E^*)$ be a VZ function with respect to a given special norm on $E \times E^*$.   Then $f$ is a VZ function with respect to all special norms on $E \times E^*$.\endslant
\Proof This is clear from the comments above and Theorem \EQthm(c).\qed
\defDefinition \INVARdef
\medbreak
\noindent
{\bf Definition \INVARdef.}\enspace Let Let $E$ be a nonzero Banach space and $f \in \PC(E \times E^*)$.   We say that $f$ is a \slant VZ function on $E \times E^*$\endslant\ if $f$ is a VZ function with respect to any one special norm on $E \times E^*$ or, equivalently, with respect to all special norms on $E \times E^*$.   This is also equivalent to the statement that $f$ is an MAS function.    
\medbreak
Theorem \BRthm(a) was obtained in \cite\VZ, Theorem 8\endcite\ under the VZ hypothesis and, in \cite\ASTWO, Theorem 4.2(2)\endcite\ under the MAS hypothesis.
\par
Theorem \BRthm(c) extends the result proved in \cite\VZ, Corollary 25\endcite\ that $\M f$ is of type (ANA).
\par
Theorem \BRthm(d) extends the result proved in \cite\ASTWO, Theorem 4.2(2)\endcite.
\par
Theorem \BRthm(f) was obtained in \cite\VZ, Corollary 7\endcite.   This is a very significant result, because maximally monotone sets $A$ of $E \times E^*$ are known such that $\overline{\pi_{E^*}(A)}$ is not convex.   \big(The first such example was given by Gossez in \cite\GOSSEZC, Proposition, p. 360\endcite\big).   Thus \big(as was first observed in \cite\VZ\endcite\big) Theorem \BRthm(f) implies that there exist maximally monotone sets $A$ that are not of the form $\M f$ for any lower semicontinuous VZ function on $E \times E^*$ or, equivalently, not of the form $\M f$ for any lower semicontinuous MAS function on $E \times E^*$.   Theorem \BRthm(f) can also be proved directly from Lemma \EXlem(a) rather than from the more circuitous argument given here.
\par
The techniques used in Theorem \BRthm\ originated in the \slant negative alignment\endslant\ analysis of \cite\BR, Section 8, pp.\ 274--280\endcite\ and \cite\HBM, Section 42, pp.\ 161--167\endcite.
\defTheorem \BRthm
\medbreak
\noindent
{\bf Theorem \BRthm.}\enspace\slant Let $E$ be a nonzero Banach space and $f \in \PCLSC(E \times E^*)$.   Assume either that $f$ is a VZ function on $E \times E^*$ or, equivalently, that $f$ is an MAS function.   Then:
\par\noindent
{\rm(a)}\enspace $\M f$ is a maximally monotone subset of $E \times E^*$.
\smallskip\noindent
{\rm(b)}\enspace Let $(x,x^*) \in E \times E^*$ and $\alpha,\beta > 0$.   Then there exists a unique value of $\omega \ge 0$ for which there exists a bounded sequence $\big\{(y_n,y_n^*)\big\}_{n \ge 1}$ of elements of $\M f$ such that,
$$\lim_{n \to \infty}\|y_n - x\| = \alpha\omega,\quad \lim_{n \to \infty}\|y_n^* - x^*\| = \beta\omega \quand \lim_{n \to \infty}\bra{y_n - x}{y_n^* - x^*} = - \alpha\beta\omega^2.$$
{\rm(c)}\enspace Let $(x,x^*) \in E \times E^* \setminus \M f$ and $\alpha,\beta > 0$.   Then there exists a bounded sequence $\big\{(y_n,y_n^*)\big\}_{n \ge 1}$ of elements of\quad $\M f \cap \big[(E \setminus \{x\}) \times (E^* \setminus \{x^*\})\big]$\quad such that,
$$\lim_{n \to \infty}\f{\|y_n - x\|}{\|y_n^* - x^*\|} = \f\alpha\beta \quand \lim_{n \to \infty}\f{\bra{y_n - x}{y_n^* - x^*}}{\|y_n - x\|\|y_n^* - x^*\|} = -1.\meqno\NAPAIRone$$
In particular, $\M f$ is of type (ANA) {\rm\big(see \cite\HBM, Definition 36.11, p.\ 152\endcite\big)}.
\smallskip\noindent
{\rm(d)}\enspace Let $(x,x^*) \in E \times E^* \setminus \M f$, $\alpha,\beta > 0$ and $\inf_{(y,y^*) \in \M f}\bra{y - x}{y^* - x^*} > -\alpha\beta$.  Then there exists a bounded sequence $\big\{(y_n,y_n^*)\big\}_{n \ge 1}$ in\quad $\M f \cap \big[(E \setminus \{x\}) \times (E^* \setminus \{x^*\})\big]$\quad  such that {\rm(\NAPAIRone)} is satisfied, $\lim_{n \to \infty}\|y_n - x\| < \alpha$ and $\lim_{n \to \infty}\|y_n^* - x^*\| < \beta$.   In particular, $\M f$ is of type (BR) {\rm\big(see \cite\HBM, Definition 36.13, p.\ 153\endcite\big)}.
\smallskip\noindent
{\rm(e)}\enspace Let $(x,x^*) \in E \times E^* \setminus \M f$, $\alpha,\beta > 0$ and $f(x,x^*) < \bra{x}{x^*} + \alpha\beta$.   Then there exists a bounded sequence $\big\{(y_n,y_n^*)\big\}_{n \ge 1}$ of elements of \quad $\M f \cap \big[(E \setminus \{x\}) \times (E^* \setminus \{x^*\})\big]$\quad  such that {\rm(\NAPAIRone)} is satisfied, $\lim_{n \to \infty}\|y_n - x\| < \alpha$ and $\lim_{n \to \infty}\|y_n^* - x^*\| < \beta$.
\smallskip\noindent
{\rm(f)}\enspace We define the {\rm projection maps} $\pi_E\colon E \times E^* \to E$ and $\pi_{E^*}\colon E \times E^* \to E^*$ by $\pi_E(x,x^*) := x$ and $\pi_{E^*}(x,x^*) := x^*$.   Then\quad $\overline{\pi_E(\M f)} = \overline{\pi_E(\dom\,f)}$ and $\overline{\pi_{E^*}(\M f)} = \overline{\pi_{E^*}(\dom\,f)}$.   Consequently, the sets $\overline{\pi_E(\M f)}$ and $\overline{\pi_{E^*}(\M f)}$ are convex.
\endslant
\Proof(a) is immediate from Theorem \EXthm(c).
\par
(b)\enspace Let $\tau := \sqrt{\beta/\alpha}$ and use the norm  $\|\cdot\|_{\infty,\tau}$ on $E \times E^*$.     Lemma \EXlem(b) provides us with a bounded sequence $\big\{(y_n,y_n^*)\big\}_{n \ge 1}$ of elements of $\M f$ such that
$$\lim_{n \to \infty}\big[\beta\|y_n - x\|^2/\alpha \vee \alpha\|y_n^* - x^*\|^2/\beta + \bra{y_n - x}{y_n^* - x^*}\big] = 0.$$
By passing to an appropriate subsequence, we can and will suppose that the three limits $\rho := \lim_{n \to \infty}\|y_n - x\|$, $\sigma := \lim_{n \to \infty}\|y_n^* - x^*\|$ and $\lim_{n \to \infty}\bra{y_n - x}{y_n^* - x^*}$ all exist.   Consequently, $\beta\rho^2/\alpha \vee \alpha\sigma^2/\beta + \lim_{n \to \infty}\bra{y_n - x}{y_n^* - x^*} = 0$, from which
$$\beta\rho^2/\alpha \vee \alpha\sigma^2/\beta = -\limn_{n \to \infty}\bra{y_n - x}{y_n^* - x^*} \le \rho\sigma = \sqrt{\beta\rho^2/\alpha}\sqrt{\alpha\sigma^2/\beta}.$$
It follows easily from this that $\beta\rho^2/\alpha = \alpha\sigma^2/\beta$ and $\lim_{n \to \infty}\bra{y_n - x}{y_n^* - x^*} = -\rho\sigma$.   The first of these equalities implies that $\rho/\alpha = \sigma/\beta$.   We take $\omega := \rho/\alpha = \sigma/\beta$, and it is immediate that $\omega$ has the required properties.   The uniqueness of $\omega$ was established in \cite\BR, Theorem 8.4(b), p.\ 276\endcite\ and \cite\HBM, Theorem 42.2(b), pp.\ 163--164\endcite.
\par
(c)\enspace Following on from (b), if $\omega = 0$ then $(\rho,\sigma) = (0,0)$, that is to say $\lim_{n \to \infty}y_n = x$ in $E$ and $\lim_{n \to \infty}y_n^* = x^*$ in $E^*$.   Since $\M f$ is closed, this would contradict the hypothesis that $(x,x^*) \not\in \M f$.   Thus $\omega > 0$, from which $\rho > 0$ and $\sigma > 0$.   (c) now follows by truncating the sequences so that, for all $n$, $\|y_n - x\| > 0$ and $\|y_n^* - x^*\| > 0$.
\par
(d)\enspace Continuing with the notation of (c), we have
$$-\alpha\beta < \infn_{(y,y^*) \in \M f}\bra{y - x}{y^* - x^*} \le\limn_{n \to \infty}\bra{y_n - x}{y_n^* - x^*} = -\rho\sigma,$$
from which $(\rho/\alpha)(\sigma/\beta) < 1$.   Since $\rho/\alpha = \sigma/\beta$, in fact  $\rho/\alpha < 1$. and  $\sigma/\beta < 1$, that is to say $\rho = \lim_{n \to \infty}\|y_n - x\| < \alpha$ and $\sigma = \lim_{n \to \infty}\|y_n^* - x^*\| < \beta$.   This gives (d).
\smallskip
(e) is immediate from (d) and the comment in Remark \FPHIrem\ that, for all $(x,x^*) \in E \times E^*$, $-\infn_{(y,y^*) \in \M f}\bra{y - x}{y^* - x^*} \le f(x,x^*) - \bra{x}{x^*}$.
.
\smallskip
(f)\enspace If $x \in \pi_E(\dom\,f)$ then there exists $x^* \in E^*$ such that $f(x,x^*) < \infty$, and so it follows from (e) that there exists $(y,y^*) \in \M f$ such that $\|y - x\| < 1/n$.   Consequently, $x \in \overline{\pi_E(\M f)}$.   Thus we have proved that $\pi_E(\dom\,f) \subset \overline{\pi_E(\M f)}$.   On the other hand, $\M f \subset \dom\,f$, and so $\overline{\pi_E(\M f)} = \overline{\pi_E(\dom\,f)}$.   We can prove in an exactly similar way that $\overline{\pi_{E^*}(\M f)} = \overline{\pi_{E^*}(\dom\,f)}$.   The convexity of the sets $\overline{\pi_E(\M f)}$ and $\overline{\pi_{E^*}(\M f)}$ now follows immediately.\qed
\smallskip
\defRemark \TWOrem
\medbreak
\noindent
{\bf Remark \TWOrem.}\enspace If we combine Theorem \EXthm(a) (using the norm $\|\cdot\|_{2,1}$ on $E \times E^*$) with the comments made in the proof of Theorem \BRthm(e)  we obtain the following result:  \slant Let $E$ be a nonzero Banach space, $f \in \PCLSC(E \times E^*)$, and $f$ be a VZ function on $E \times E^*$.
Then, for all $(x,x^*) \in E \times E^*$,
$$\eqalign{
\infn_{(y,y^*) \in \M f}\sqrt{\|y - x\|^2 + \|y^* - x^*\|^2}
&\le \rttwo\sqrt{-\infn_{(y,y^*) \in \M f}\bra{y - x}{y^* - x^*}}\cr
&\le \rttwo\sqrt{f(x,x^*) - \bra{x}{x^*}}.}$$\endslant
This strengthens the result proved in \cite\VZ, Theorem 4\endcite, namely that
$$\infn_{(y,y^*) \in \M f}\sqrt{\|y - x\|^2 + \|y^* - x^*\|^2} \le 2\sqrt{f(x,x^*) - \bra{x}{x^*}}.$$
As we observed in Remark \FPHIrem, the constant $\sqrt 2$ is sharp.
\defDefinition \SRdef
\medbreak
\noindent
{\bf Definition \SRdef.}\enspace Let $E$ be a nonzero Banach space and $A$ be a nonempty monotone subset of  $E \times E^*$.   We say that $A$ is \slant of type (NI)\endslant\ if, for all $(x\dbs,x^*) \in E\dbs \times E^*$,
$$\infn_{(s,s^*) \in A}\bra{x^* - s^*}{x\dbs - \wh s} \le 0.$$
This concept was introduced in \cite\RANGE, Definition 10, p.\ 183\endcite. 
We say that $A$ is \slant strongly representable\endslant\ if there exists $f \in \PCLSC(E \times E^*)$ such that\quad $f \ge q$ on $E \times E^*$, \quad $f^* \ge \qt$ on $E\dbs \times E^*$\quad (i.e., $f$ is a lower semicontinuous MAS function)\quad and $\M f = A$.   This concept was introduced and studied in \cite\ASTWO\endcite, \cite\ASTHREE\endcite\ and \cite\VZ\endcite. 
\medbreak
Theorem \SRthm\ was motivated by and extends that proved in \cite\ASTHREE, Theorem 1.2\endcite.   The most significant part of it is the fact that (a) implies (d) and (a) implies (e).   In particular, if $A$ is maximally monotone of type (NI), then the conclusions of Theorem \BRthm(b--f) hold (with $\M f$ replaced by $A$).   This leads to a substantial generalization of \cite\BR, Theorem 8.6, pp.\ 277--278\endcite\ and \cite\HBM, Theorem 42.6, pp.\ 163--164\endcite.   The fact that $\overline{\pi_EA}$ and $\overline{\pi_{E^*}A}$ are convex whenever $A$ is of type (NI) was first proved by Zagrodny in \cite\ZAGRODNY\endcite.
\defTheorem \SRthm
\medbreak
\noindent
{\bf Theorem \SRthm.}\enspace\slant Let $E$ be a nonzero Banach space and $A$ be a maximally monotone subset of $E \times E^*$.   Then the following conditions are equivalent:
\smallskip\noindent
{\rm(a)}\enspace $A$ is of type (NI).
\smallskip\noindent
{\rm(b)}\enspace For all $(x\dbs,x^*) \in E\dbs \times E^*$,\quad $\supn_{(s,s^*) \in A}\big[\bra{s}{x^*} + \bra{s^*}{x\dbs} - \bra{s}{s^*}\big] \ge \bra{x^*}{x\dbs}$.
\smallskip\noindent
{\rm(c)}\enspace For all $(x\dbs,x^*) \in E\dbs \times E^*$,\enspace $\supn_{(y,y^*) \in E \times E^*}\big[\bra{y}{x^*} + \bra{y^*}{x\dbs} - \Phi_A(y,y^*)\big] \ge \bra{x^*}{x\dbs}$.
\smallskip\noindent
{\rm(d)}\enspace $A$ is strongly representable.
\smallskip\noindent
{\rm(e)}\enspace There exists a lower semicontinuous VZ function on $E \times E^*$ such that $\M f = A$.
\smallskip\noindent
{\rm(f)}\enspace $\Phi_A$ is a VZ function on $E \times E^*$.
\smallskip\noindent
{\rm(g)}\enspace $^*\Theta_A$ is a VZ function on $E \times E^*$.
\smallskip\noindent
{\rm(b$_1$)}\enspace If $h \in \PC(E \times E^*)$ and\quad $^*\Theta_A \ge h$ on $E \times E^*$\quad then, for all $(x\dbs,x^*) \in E\dbs \times E^*$,
$$h^*(x\dbs,x^*) = \supn_{(y,y^*) \in E \times E^*}\big[\bra{y}{x^*} + \bra{y^*}{x\dbs} - h(y,y^*)\big] \ge \bra{x^*}{x\dbs}.\meqno\SRone$$
\smallskip\noindent
{\rm(b$_2$)}\enspace If $h \in \PCLSC(E \times E^*)$ and\quad $^*\Theta_A \ge h \ge \Phi_A$ on $E \times E^*$\quad then, for all $(x\dbs,x^*) \in E\dbs \times E^*$, {\rm(\SRone)} is satisfied.
\smallskip\noindent
{\rm(c$_1$)}\enspace There exists $h \in \PCLSC(E \times E^*)$ such that\quad $^*\Theta_A \ge h \ge \Phi_A$ on $E \times E^*$\quad and, for all $(x\dbs,x^*) \in E\dbs \times E^*$, {\rm(\SRone)} is satisfied.
\smallskip\noindent
{\rm(c$_2$)}\enspace There exists $h \in \PC(E \times E^*)$ such that \quad $h \ge \Phi_A$ on $E \times E^*$\quad and, for all $(x\dbs,x^*) \in E\dbs \times E^*$, {\rm(\SRone)} is satisfied.
\endslant
\Proof These results are all immediate from the corresponding parts of Theorem \MULTIthm.\qed 
\defSection \FMsec
\bigbreak   
\centerline{\bf \FMsec\quad Appendix: a nonhausdorff Fenchel--Moreau theorem}
\medskip\noindent
In Remark \THrem, we referred to the Fenchel--Moreau theorem for (possibly nonhausdorff) locally convex spaces.   We shall give a proof of this result in Theorem \FMthm.   When we say that $X$ is a \slant locally convex space\endslant, we mean that $X$ is a nonzero real vector space endowed with a topology compatible with its vector structure with  a base of neighborhoods of $0$ of the form $\big\{x \in X\colon\ S(x) \le 1\big\}_{S \in \Sem(X)}$, where $\Sem(X)$ is a family of seminorms on $X$ such that if $S_1 \in \Sem(X)$ and $S_2 \in \Sem(X)$ then $S_1 \vee S_2 \in \Sem(X)$; and if $S \in \Sem(X)$ and $\lambda \ge 0$ then $\lambda S \in \Sem(X)$.   If $L$ is a linear functional on $X$ then $L$ is continuous if, and only if, there exists $S \in \Sem(X)$ such that $L \le S$ on $X$.
\par
As an example of the construction above, we can suppose that $X$ and $Y$ are vector spaces paired by a bilinear form $\bra\cdot\cdot$.   Then $\big(X,w(X,Y)\big)$ is a locally convex space with determining family of seminorms $\big\{|\bra{\cdot}{y_1}| \vee \cdots \vee |\bra{\cdot}{y_n}|\big\}_{n \ge 1,\ y_1,\dots, y_n \in Y}$.
\par
The author is grateful to Constantin Z\u alinescu for showing him a proof of Theorem \FMthm\ based on the standard (Hausdorff) result and a quotient construction.   The proof we give here is a simplification of the result on ``Fenchel--Moreau points'' of \cite\HBL, Theorem 5.3, pp.\ 157--158\endcite\ or \cite\HBM, Theorem 12.2, pp.\ 59--60\endcite, which is also valid in the nonhausdorff setting.   
\defTheorem \FMthm
\medbreak
\noindent
{\bf Theorem \FMthm.}\enspace\slant Let $X$ be a locally convex space and $f \in \PC(X)$ be lower semicontinuous.   Write $X^*$ for the set of continuous linear functionals on $X$.   If $L \in X^*$, define $f^*(L) := \sup_X\big[L - f\big]$.    Let $y \in X$.  Then
$$f(y) = \supn_{L \in X^*}\big[L(y) - f^*(L)\big].\meqno\FMone$$\endslant
\Proof Since, for all $L \in X^*$,
$L(y) - f^*(L) = \infn_{x \in X}\big[L(y) - L(x) + f(x)\big] = (f \episum L)(y)$ and the inequality ``$\ge$'' in (\FMone) is obvious from the definition of $f^*(L)$, we only have to prove that
$$f(y) \le \supn_{L \in X^*}(f \episum L)(y)\big].\meqno\FMtwo$$   
Let $\lambda \in \r$  and $\lambda < f(y)$.   Since $f$ is proper, there exists $z \in \dom\,f$.  Choose $Q \in \Sem(X)$ such that
$$Q(z - x) \le 1 \qlr f(x) > f(z) - 1\meqno\FENMORSUFFone$$
and
$$Q(y - x) \le 1 \qlr f(x) > \lambda.\meqno\FENMORSUFFtwo$$
\par
We first prove that
$$(f \episum Q)(z) \ge f(z) - 1.\meqno\FENMORSUFFthree$$
To this end, let $x$ be an arbitrary element of $X$.  If $Q(z - x) \le 1$ then (\FENMORSUFFone) implies that $f(x) + Q(z - x) \ge f(x) > f(z) - 1$.   If, on the other hand, $Q(z - x) > 1$, let $\gamma := 1/Q(z - x) \in \,]0,1[\,$ and put $u := \gamma x + (1 - \gamma)z$.   Then $Q(z - u) = \gamma Q(z - x) = 1$ and so, from the convexity of $f$, and (\FENMORSUFFone) with $x$ replaced by $u$,
$$\gamma f(x) + (1 - \gamma)f(z) \ge f\big(\gamma x + (1 - \gamma)z\big) = f(u) > f(z) - 1.$$
Substituting in the formula for $\gamma$ and clearing of fractions yields $f(x) + Q(z - x) \ge f(z)$.   This completes the proof of (\FENMORSUFFthree).
\par
Now let $M \ge 1$ and $M \ge \lambda + 2 + Q(z - y) - f(z)$.   We will prove that   
$$(f \episum MQ)(y) \ge \lambda.\meqno\FENMORSUFFfour$$
To this end, let $x$ be an arbitrary element of $X$.  If $Q(y - x) \le 1$ then (\FENMORSUFFtwo) implies that $f(x) + MQ(y - x) \ge f(x) > \lambda$.   If, on the other hand, $Q(y - x) > 1$ then, from (\FENMORSUFFthree),
$$\eqalign{f(x) + MQ(y - x)
&= f(x) + Q(y - x) + (M - 1)Q(y - x)\cr
&\ge f(x) + Q(z - x) - Q(z - y) + (M - 1)\cr
&\ge f(z) - 1 - Q(z - y) + M - 1 \ge \lambda,}$$
which completes the proof of (\FENMORSUFFfour).   The ``Hahn--Banach--Lagrange theorem'' of\break \cite\HBL, Theorem 2.9, p.\ 153\endcite\ or \cite\HBM, Theorem 1.11, p.\ 21\endcite\ now provides us with a linear functional $L$ on $X$ such that $L \le MQ$ on $X$ and $(f \episum L)(y) \ge \lambda$.   (\FMtwo) now follows by letting $\lambda \to f(y)$.\qed
\bigskip
\centerline{\bf References}
\nmbr\BS
\item{[\BS]} R. S. Burachik and B. F. Svaiter, \slant Maximal monotonicity, conjugation and the duality product\endslant,  Proc. Amer. Math. Soc. {\bf 131}  (2003), 2379--2383.
\nmbr\FITZ
\item{[\FITZ]} S. Fitzpatrick, \slant Representing monotone operators
by convex functions\endslant,   Workshop/ Miniconference on
Functional Analysis and Optimization (Canberra, 1988),  59--65, Proc.
Centre Math. Anal. Austral. Nat. Univ., {\bf 20}, Austral. Nat. Univ.,
Canberra, 1988.
\nmbr\GOSSEZC
\item{[\GOSSEZC]}J.- P. Gossez, \slant On a convexity property of the range of a maximal monotone operator\endslant, Proc. Amer. Math. Soc. {\bf 55} (1976), 359--360.
\nmbr\ASTWO
\item{[\ASTWO]}M. Marques Alves and B. F. Svaiter, \slant Br\o ndsted--Rockafellar property and maximality of monotone operators representable by convex functions in non--reflexive Banach spaces.\endslant, http://arxiv.org/abs/0802.1895v1, posted Feb 13, 2008, to appear in J. of Convex Anal.
\nmbr\ASTHREE
\item{[\ASTHREE]}M. Marques Alves and B. F. Svaiter, \slant A new old class of maximal monotone operators.\endslant, http://arxiv.org/abs/0805.4597v1, posted May 29, 2008.
\nmbr\MT
\item{[\MT]}J.--E. Mart\'\i nez-Legaz and M. Th\'era, \slant $\eps$--Subdifferentials in terms of subdifferentials\endslant, Set--Valued
Anal. {\bf 4} (1996), 327--332.
\nmbr\PENOT
\item{[\PENOT]}J.--P. Penot, \slant The relevance of convex analysis for the study of monotonicity\endslant,  Nonlinear Anal. {\bf 58}  (2004), 855--871.
\nmbr\FENCHEL
\item{[\FENCHEL]}R. T. Rockafellar, \slant Extension of Fenchel's duality theorem for convex functions\endslant, Duke Math. J. {\bf33} (1966), 81--89.
\nmbr\RANGE
\item{[\RANGE]}S. Simons, \slant The range of a monotone operator\endslant, J. Math. Anal. Appl. {\bf 199} (1996), 176--201.
\nmbr\BR
\item{[\BR]}-----, \slant Maximal monotone multifunctions of Br{\o}ndsted--Rockafellar type\endslant, Set--Valued Anal. {\bf7} (1999),
255--294.
\nmbr\PANDM
\item{[\PANDM]}-----, \slant Positive sets and Monotone sets\endslant, J. of Convex Anal., {\bf 14} (2007), 297--317.
\nmbr\HBL
\item{[\PANDM]}-----, \slant The Hahn--Banach--Lagrange theorem\endslant, Optimization, {\bf 56} (2007), 149--169.
\nmbr\HBM
\item{[\HBM]}-----, \slant From Hahn--Banach to monotonicity\endslant, 
Lecture Notes in Mathematics, {\bf 1693},\break second edition, (2008), Springer--Verlag.
\nmbr\ZAGRODNY
\item{[\ZAGRODNY]}D. Zagrodny, \slant The convexity of the closure of the domain and the range of a maximal monotone multifunction of Type NI\endslant, to appear in Set--Valued Anal.
\nmbr\VZ
\item{[\VZ]}M. D. Voisei and C. Z\u{a}linescu, \slant Strongly--representable operators\endslant, http://arxiv.org/ abs/0802.3640v1, posted February 25, 2008. 
\nmbr\ZBOOK
\item{[\ZBOOK]} C. Z\u{a}linescu, \slant Convex analysis in general vector spaces\endslant, (2002), World Scientific.
%
%
\Signoff
\bigskip
Department of Mathematics\par
University of California\par
Santa Barbara\par
CA 93106-3080\par
U. S. A.\par
email:  simons@math.ucsb.edu
\bye